\documentclass[11pt,english]{article}

\usepackage{graphicx,psfrag,epsfig,amsfonts,amssymb,amsmath,babel}

 \usepackage[colorlinks=true]{hyperref}

\hypersetup{urlcolor=blue, citecolor=red}

\newtheorem{theorem}{Theorem}[section]
\newtheorem{corollary}[theorem]{Corollary}

\newtheorem{lemma}[theorem]{Lemma}
\newtheorem{proposition}[theorem]{Proposition}

\newtheorem{definition}[theorem]{Definition}
\newtheorem{remark}[theorem]{Remark}

\begin{document}

\title{Dynamics of  large   cooperative pulsed-coupled networks}

\author{Eleonora Catsigeras\thanks{Instituto de Matem\'{a}tica y Estad\'{\i}stica Rafael Laguardia (IMERL),
 Fac. Ingenier\'{\i}a,  Universidad de la Rep\'{u}blica,  Uruguay.
 E-mail: eleonora@fing.edu.uy
  EC was partially supported by CSIC of Universidad de la Rep\'{u}blica and ANII of Uruguay.} }

\date{}

\maketitle

\begin{abstract}
We study the deterministic dynamics of   networks ${\mathcal N}$ composed by $m$  non identical, mutually pulse-coupled cells. We assume
weighted, asymmetric and positive (cooperative) interactions among the cells, and arbitrarily large values of $m$. We consider two cases of  the network's graph: the complete graph, and the existence of a large    core (i.e. a large complete subgraph). First, we prove that  the system periodically  eventually synchronizes  with a natural \lq\lq spiking period\rq\rq \ $p \geq 1$, and that if the cells are mutually structurally identical or similar, then the synchronization is complete ($p= 1$) .   Second, we prove that the   amount of information $H$  that ${\mathcal N}$   generates or processes, equals   $\log p$. Therefore, if ${\mathcal N}$ completely synchronizes, the information   is null. Finally, we prove that ${\mathcal N}$  protects the cells from their risk of death.
\end{abstract}

{\noindent \footnotesize {{\em MSC }2010:  {\  Primary: 37NXX, 92B20; Secondary: 34D06, 05C82, 94A17, 92B25} }}

\section{Introduction}

The  theory of deterministic dynamical systems composed by two or more  coupled dynamical units  assumes that each unit - which we call \em cell \em- has a proper own dynamics, and that the   couplings among the units are \em interactions \em  that depend  on the instantaneous states of the cells \cite{Young2008,ChazottesFernandez2005}.

Among the systems of interacting units, we focus on those that are pulsed-coupled (i.e. the interactions are instantaneous). In particular, the  global system   can be understood as a  multi-dimensional differential or difference equation with impulsive terms  \cite{MassBishop2001}.

On the one hand, the theory of dynamically interacting units is a source of   mathematical  open questions     \cite{Young2008}. In particular, those systems governed by impulsive differential equations    \cite{stamov2007}  pose mathematical problems that are mostly open, except in particular cases or low dimensions.

On the other hand, the network of interacting dynamical units is a particular model of a    coalitional game, that evolves or changes on time. We will discuss this relation in Subsection \ref{subsectionNew_CooperativeGames}.

 The dynamical systems composed by mutually interacting units, and in particular those that are pulsed-coupled, have relevance in many applications.  As examples  in Physics, a two-dimensional impulsive differential equation  models the joint dynamics of two or more coupled oscillators  \cite{CabezaMartiRubidoKahan2011}. In particular, they are used in applications to  Light-Controlled-Oscillators (LCO),  \cite{LCOrecomMartiCabeza2_2003}. The mathematical investigation of a genetic regulatory network of two antagonist genes, is also modeled as a network of   pulse-coupled units \cite{arnaud}.
In Neuroscience,  among the theoretical methods of research,   the mathematical analysis of the dynamics of pulse-coupled networks is applied  \cite{ErmentroutTerman2010, Izhikevich2007}.
In Engineering,   networks of coupled dynamical units are designed for     control systems and communications  \cite{YangChua}. Computational research on artificial intelligence, by means of   artificial neuronal networks, is used to analyze, simulate, and investigate on data obtained from dynamical systems of interacting units with a large degree of complexity  \cite{cottrell2012}.
  In Economics,   networks of coupled units are used to investigate the equilibrium states in social systems of  interacting agents \cite{elvio2009}. Artificial neuronal networks are applied for the prediction  of the exchange market  \cite{misas},   also to investigate on financial markets  \cite{finantialMarkets}, and for the accounting of financial applications \cite{TesisFinlandia}.
  In Ecology, networks of  interacting units model the dynamics of predator-prey communities of two or more species \cite{ecology2009,ecologiaMutualInterference3_2010}.
  The dynamics of  an infectious disease, taking into account the interaction among populations of diverse infectious agents, is   modelled as a  neural network  \cite{infectiousdiseasemodel2010}.
  In Geosciences, the forecast of ozone peaks in weather prediction uses computational methods on artificial neuronal networks  \cite{OzonePeaks}. In Social Sciences, the dynamics of large WWW social networks is mathematically modelled by the interactions of their individuals \cite{wattsStrogatzNewman2002,wattsStrogatz1998}; and the self- synchronization of many small clusters of cells in a low-dimensional network models the dynamics of a medieval social network \cite{socialNetworks}.

\subsection{The object and method of research} \label{subsection1.1}

Along this paper we investigate, by exact mathematical analysis and deductive proofs, the global dynamics of certain   pulse-coupled deterministic networks of $m $ cooperative   cells, for   large values of $m$.

Each cell $i \in \{1, 2, \ldots, m\}$ is   governed by a deterministic dynamical sub-system, which - if $i$ were hypothetically isolated from the network - we call \em the free dynamics \em of $i$. Besides, each cell $i$  \em acts on \em the other cells $j \neq i$ of the network  at certain instants $t_i$, which we call the \em spiking   or   milestone instants \em of  $i$. Conversely, each cell $i$   \em receives the actions from \em the other cells $j \neq i$ at the spiking instants $t_j$ of $j$.

 The free dynamics of each cell   evolves governed by a finite-dimensional ordinary \em differential equation, \em joint with an autonomous and instantaneous \em reset or update rule. \em The update   rule   applies when the state of $i$ arrives to a pursued \em goal or threshold level $\theta_i$. \em The update rule   resets the state of $i$, or equivalently, it   changes the velocity according to which the free dynamics evolves.   The  free dynamics of the many cells of the network may be mutually    different. As a particular case, in Neuroscience  the model of the free dynamics of each cell - by integration of a differential equation plus a reset or update rule - is called \em integrate and fire. \em
 Each reset or update event of the cell (the neuron) is called \em a spike.  \em    In brief, the spike of a cell $i$ is produced when its state arrives to the goal   $\theta_i$.

The cells compose the network  by \em mutual interactions \em between any ordered pair $(i,j)$ such that   $1 \leq i,j \leq m$ and $i \neq j$. These interactions    exist in both directions (some   interactions may be zero), but are neither necessarily symmetric nor simultaneous. Roughly speaking,
the     action from the cell $i$ to the cell $j \neq i$ is a discontinuity jump $\Delta_{i,j}$ - applied on   the state of $j$ - that is produced at the instant $t_i$  when $i$ spikes. So, the instant $t_i$ depends on the state of $i$. The rule is the same to define the action from the cell $j$ to the cell $i$, but the matrix $(\Delta)_{i,j}$ is not necessarily symmetric. Besides,   the instant $t_i$ when the action $\Delta_{i,j}$ is applied, is in general different from the instant $t_j$ when the action $\Delta_{j,i}$ is applied.

\vspace{.3cm}

Our purpose of research is to find qualitative and quantitative relevant characteristics of the global dynamics of such an abstract network, while time $t$ evolves to the future. As said above, the methodology is to find the exact abstract mathematical statements and their   deductive proofs.

\vspace{.3cm}

We take   the main ideas from \cite{YoPierre}, in which  a model of a network,  composed by integrate and fire biological neurons,  is studied by the exact mathematical method.

The main differences between the dynamical system that we study here and the one studied in \cite{YoPierre}, are the following:
 First, we assume that     the cells are cooperative (which in Neuroscience are called excitatory). This means that   $\Delta_{i,j} \geq 0$ for all ordered pairs $(i,j)$ such that $i \neq j$, and the value zero is admitted in some of our results. In \cite{YoPierre}, any sign of $\Delta_{i,j}$ is admitted by hypothesis, but only nonzero values are assumed. Second, we do not assume that the free dynamics   is the same for all the cells. In \cite{YoPierre} all the cells are identical.   Third, we neither assume the linearity of the differential equation  that governs  the free dynamics of each cell,  nor the existence of a   Lyapunov stable   equilibrium state for the solution flow of this differential equation. In \cite{YoPierre} these latter two conditions are assumed.

\subsection{The network as a cooperative game that evolves on time} \label{subsectionNew_CooperativeGames} The   mathematical model   of the network   ${\mathcal N} $ that we study   in this paper is also    an evolutive game  represented by a graph whose $m$ vertices are the players $i \in \{1, \ldots, m\}$ (the cells), and whose edges $\Delta_{i,j}$ (the interactions) are directed and weighted. The hypothesis of cooperativity among the individuals, namely $\Delta_{i,j} \geq 0$ for all $i \neq j$ makes the game work in an \em imitate the best \em strategy, which produce players that  adopt a myopic behaviour (\cite{golman}). In fact, by hypothesis, each cell or player  $i$  just knows its own actions to the other players $j \neq i$, the value of its own satisfaction variable, and the actions it receives from the other cells. But $i$ ignores   the global state and dynamical behaviour of the whole network.

The model of the network ${\mathcal N}$, as described in Subsection \ref{subsection1.1} from the dynamical viewpoint, is a cooperative or coalitional game, that changes on time. It disregards the individual strategies of its cells (or players) and instead, it focusses on the coalitions (which we call clusters of cells), defined as nonempty subsets of cells $i \in \{1, \ldots, m\}$ whose satisfaction variables $S_i$ arrive simultaneously to their respective goals or threshold levels $\theta_i >0$.

The characteristic function $\nu: 2^m \mapsto \mathbb{R}^+$ of the coalitional game (i.e. the function  assigning the total gain or payment $\nu(A)$ to each   coalition $A \subset \{1, \ldots, m\}$, with the agreement $\nu(\emptyset)= 0$) can be defined as the sum of the goal levels $\theta_i$ of the cells $i \in A$. In fact, at each instant $t_n$ for which a nonempty coalition $I_n$ is formed (i.e. a cluster of simultaneously spiking cells is exhibited), the satisfaction variable  of each cells   $i \in I_n$ equals its respective goal or threshold level $\theta_i$.  Since $\theta_i >0$ for all $i \in {\mathcal N}$,   the coalitional game is \em convex\em, i.e.
\begin{equation}
\label{eqnConvexity} \nu(A \cup B) + \nu(A \cap B) \geq \nu(A) + \nu(B) \ \ \forall \ A, B \subset {\mathcal N}
\end{equation}
We note that, since in our case $\nu(A)= \sum_{i \in {\mathcal A}} \theta_i$, inequality (\ref{eqnConvexity}) is indeed an equality.
Thus:
$$\nu(A \cup \{i\}) - \nu(A)  \leq \nu(B \cup \{i\}) - \nu(B) \ \ \forall \ A \subset B \subset {\mathcal N}\setminus \{i\}, \ \ \forall \ i \in {\mathcal N}.$$
In other words, the gain for a player $i$ for joining a coalition $B$ larger than $A$, is non negative.

As for any convex coalitional game, ${\mathcal N}$ has a nonempty core of solutions. In Game Theory, a solution in the core is a vector $$(z_1, z_2, \ldots, z_i, \ldots, z_m) \in \mathbb{R}^m,$$ which is called \lq\lq allocation\rq\rq \ or \lq\lq payoff vector\rq\rq, such that
\begin{equation}\label{eqnCoreGameTheory} \sum_{i \in I} z_i \geq \nu(I) \ \forall \ I \subset {\mathcal N}.\end{equation}

For our model, the payoff vector at the $n$-th. instant $t_n$ when at least one cell arrives to its goal level and spikes, can be defined by the following formula:
$$z_i = S_i(t_n^-) + \sum_{j \neq i, \ j \in I_n} \Delta_{j,i},$$
where $S_{i}(t_n^-) \in [0, \theta_i]$ is the value of the satisfaction variable $S_i$ of the cell $i$ just before instant $t_n$ (cf. Definition \ref{definitionSpikeMilestoneSatisfactionVariable}),
$I_n$ is the cluster at the spiking instant $t_n$ (cf. Definition \ref{definitionCluster}) and $\Delta_{j,i} \geq 0$ is the cooperative action from the cell $j \in I_n$ to $i \neq j$ at any instant for which $j$ spikes. In Equation (\ref{eqnInteraction}) we will precisely state the rule according to which the cooperative interactions $\Delta_{j,i}$ among the cells of the network increase their satisfaction variables. In brief, the cell $i$ spikes at instant $t_n$, hence it belongs to the cluster $I_n$, if and only if one of the following conditions is satisfied: either the satisfaction variable $S_i$ arrives to the goal value $\theta_i$ spontaneously (due to the free dynamics of $i$) at instant $t_n^-$, or the satisfaction variable at instant $t_n^-$ is smaller than $\theta_i$  but   suddenly increases to become larger or equal than $\theta_i$, by the adding of some positive cooperative actions $\Delta_{j,i}$ from cells $j \neq i$ that spike at the same instant $t_n$. In any case, the component $z_i$ of the payoff vector is $z_i \geq \theta_i$  if $z_i \in I_n$. Since $\nu(I_n) = \sum_{i \in I_n} \theta_i$ we deduce that Inequality (\ref{eqnCoreGameTheory}) holds in the coalitional game ${\mathcal N}$.

Although   the characteristic function $\nu$ is assumed to be invariant with time $t$, the payoff vector changes with $t$, because the satisfaction variables evolve on time, and the interactions among the cells are not constantly applied at any time. In fact $\Delta_{j,i}$ is effectively added to the satisfaction variable $S_i$ just at the instants for which the cell $j$ spikes. So ${\mathcal N}$ is an evolutive coalitional game: its solutions depend on time.

In Definition \ref{definitionSincronizacion} we define the global synchronization phenomenon of the network ${\mathcal N}$ as the recurrent exhibition - infinitely many times in the future - of   the so called \lq\lq grand coalition\rq\rq.  Precisely,  all the cells of the network compose a single cluster: all of them spike simultaneously, and this phenomenon occurs infinitely many times. In Theorem \ref{theoremFullCooperativism} we prove that, if the number $m$ of players is large enough in relation with the minimum positive interaction (namely, $m$ is large compared with the minimum payoff component at any instant), then the global synchronization occurs recurrently. Nevertheless,   between two  instants   $t^*_2 > t^*_1$ when the grand coalition is exhibited, there may appear many coalitions or clusters that are smaller than the grand coalition. In Theorem \ref{TheoremPart(e)} we prove that if more different coalitions recurrently appear in the future, then the amount of information of the global dynamics of the network increases.

  In spite that the network can be studied as a coalitional evolutive game, it has some important differences    from the    models that are usually investigated in Game Theory, even from those that focus on evolutive games (e.g. \cite{golman,elvio2009,miltaich}):   For instance, in our model,  each cell or player $i$  decides \em the instants when it   acts \em in the network, by integrating the states of   its own internal dynamics with the actions (or payoffs) that it had received from the other players. So,  there does not exist a forcing external mechanism making   the players interact   sharing a resource. As a counterpart, neither the  cell nor the coalition, has an alternative  to decide how to choose among several strategies, its actions on the network. The cells and the coalitions are not   optimizers. The network  is just a deterministic dynamical system whose state suffers instantaneous changes according to the states of its cells   and, as a result of the positive interactions among the cells, it shows   patterns of coalitions   that spike simultaneously infinitely many times in the future.

\subsection{The hypothesis} \label{subsubsectionHypothesis}

We adopt three types of hypothesis. First, we assume that the number of cells - in a part or the whole network  - is large enough with respect to certain other parameters of the network (Inequality (\ref{eqnm'>Theorem1}) of Theorem \ref{theoremFullCooperativism} and a similar inequality in Corollary \ref{corollary1}). Second, we study  the cases for which  the  cells  are   cooperative (cf. Definition \ref{definitioncooperative}). Third, we assume    that  the network's graph is  complete, or, if not, that there exists   a large \em   cooperative core \em  (i.e., a complete subgraph of   cells which are fully cooperative, cf. Definition \ref{definitionCore}).

\subsection{The results to be proved}
Under  the hypothesis above described, we prove   three   results: the periodic eventual synchronization, the positiveness of the protection factor of the network against the risk of death of its cells, and the value of the total amount of information that the network is able to generate or process.

\subsubsection{Results about synchronization}

In Theorems \ref{theoremFullCooperativism}, \ref{TheoremPart(c)}  and Corollary \ref{corollary1}, we   prove   that the  network   necessarily self-synchronizes the spikes of   clusters   of cells periodically, after a transitory waiting time, and that this behaviour is robust (i.e. it persists under small perturbation of the numerical and functional parameters of the network).  After the transients, the synchronized clusters   are independent of the  initial state. But the finite transitory time-interval itself,  does depend  on the initial state.

We give  the word \lq\lq   synchronization\rq\rq \   an exact mathematical meaning (cf. Definitions \ref{definitionSincronizacionPeriodica} and \ref{definitionSincronizacion}). By eventual \em periodic spike-synchronization, \em we mean  that  the interactions among  \em all \em the cells within the same cluster are produced simultaneously, infinitely many times in the future and periodically - but \em not necessarily at all the  spiking instants of the network \em   (Definition \ref{definitionSincronizacionPeriodica}). Nevertheless, this definition  is rather non-standard   in Physics. (For the concept  of  synchronization and phase locking in Physics see for instance \cite{Pikovsky}.) In fact, the classical definition of synchronization     requires    all the cells    be identical (e.g. \cite{MirolloStrogatz}). But along this paper, we are not assuming the hypothesis of identical cells, neither in the whole network nor in a part of it.

Joining the results of  Theorems \ref{theoremFullCooperativism} and Theorem \ref{TheoremParts(fgh)},  we prove that the events' synchronization of more than one cluster  is produced only if the cells are mutually     different. In other words, if all the cells are structurally  identical - also if they are non identical but similar - we prove that  all
the events of the whole network synchronize, after a transitory time-interval and from any initial state (Theorem \ref{TheoremParts(fgh)}, Part (a)). In brief, networks of similar cooperative cells  form a unique spike-synchronized cluster.

In their seminal article  \cite{MirolloStrogatz}, Mirollo and Strogatz proved that networks of identical fully cooperative cells, with constant interaction $\Delta_{i,j} = \Delta >0 \ \forall \ i \neq j$, synchronize all their spikes.
In \cite{Bottani1996}, Bottani proved a very general result under a certain stability hypothesis,  stating the synchronization of complete networks composed by a large number of excitatory neurons (i.e. cooperative  cells) that are modelled as integrate and fire oscillators: Either all the oscillators evolve synchronized in block, or subsets of synchronized oscillators appear always in stable avalanches.
The synchronization of  several clusters of cells, provide   repetitive patterns that allow the network to organize the information.   In Theorem  \ref{theoremFullCooperativism}     we generalize   Bottani's results to mathematically abstract networks of pulse-coupled cooperative cells, that are neither necessarily  identical nor similar, and without the stability hypothesis. In Corollary \ref{corollary1} we generalize the statements for some kind of not completely connected networks.

\subsubsection{Results  about the protection and the amount of information} In Theorem \ref{TheoremPart(d)}   we prove that cooperative networks with a large number of cells protect their cells from external negative interferences, diminishing their risk  of death (cf. \cite{ErmentroutKopell1990}). This protection is due to the synchronization of sufficiently large clusters. As a counterpart, in  Part (c) of Theorem \ref{TheoremParts(fgh)}, we prove that if a complete synchronization is achieved, then the total amount of information that the network can generate,  or process, is null.  But if the synchronization is not complete - as said above this occurs if  perceptible differences   exist among the free dynamics of the cells - then the network  is able to process   a positive amount of information  (see Theorem \ref{TheoremPart(e)}, in the case that the spiking period  $p$ is not 1). Nevertheless,   this amount of information is necessarily much lower than the theoretical maximum that a general network could process, according to the number of cells that are employed in the task. This is because, due to Theorem \ref{theoremFullCooperativism}, in large cooperative networks there exist large clusters of cells whose spikes are produced all together. Therefore, they    contribute to enlarge the number of different spiking-code patterns, as if each cluster had a unique cell.

 The null amount of information of the   cooperative networks that completely synchronize, opposes  to   many networks composed by all antagonist cells. These latter may exhibit a large amount of (rather unpredictable) information. See for instance the \lq\lq virtual chaos\rq\rq \  in \cite{Cessac, cessacVieville}, or the \lq\lq stable chaos\rq\rq \  in \cite{PolitiTorcini2010}.  As a counterpart, the networks composed by all antagonist cells  unprotect their cells and so, they enlarge  their  risk of death (cf. Remark  \ref{remarkProtectionFactor}). As a consequence, if all the interactions were antagonist, the whole network would be in risk, since the death of many of its cells can drastically diminish the initial richness of its global dynamics.

 Therefore,   a large network that performs both features: protects their cells and processes a large amount of \em not always structured and predictable \em information, should theoretically be composed by some cooperative cells and by some antagonist ones, in a more or less balanced interplay   among them.   Nevertheless,   networks of antagonist cells (which in Neuroscience are called inhibitory neurons) also  may synchronize and exhibit null amount of information,  if the hypothesis of instantaneous spikes is released to take into account the delays  \cite{CvanVreeswijk}.

 Large networks of all cooperative cells with instantaneous interactions-  the ones that we study along this paper - also show a balanced interplay between the amount of information and the protection of the network against the risk of death of its cells. In fact, in Theorem \ref{TheoremParts(fgh)}, we prove that  the network maximizes its protection to the cells, if they have mutually similar free dynamics. But if so,    the total amount of information  becomes   zero (Part (c) of Theorem \ref{TheoremParts(fgh)}). In other words, if the network is able to process a  positive amount of information, then the cells are mutually different. In particular, the time-constants of their respective free dynamics must be   diverse. In such a case, the network may show  many different clusters of mutually synchronized cells, and a long  period $p$ until the spike patterns of the whole network repeat. In Theorem \ref{TheoremPart(e)}, we prove  that the  amount of information $H$ that a fully cooperative network is able to generate or process, equals   $\log_2 p$. So, $H$ is larger if the number of synchronized clusters  and  the period $p$  are larger.

 \vspace{.3cm}

The paper is organized as follows:

In Section \ref{subsectionBasicDefinitions} we  pose the previous  definitions and hypothesis that are assumed for the mathematical model of a cooperative network of pulsed coupled dynamical units. In Section \ref{subsectionStatementOfResults}, we include   the   particular  mathematical definitions and the    statements of the   theorems that we will prove along the paper.   From Section  \ref{sectionSincronizacion} to Section \ref{sectionProofOfCorollaries}  we write the  proofs.

 \section{Previous definitions and assumptions}
\label{subsectionBasicDefinitions}

We call the scale of the cell the \em micro-scale \em  relative to the whole network ${\mathcal N}$ under study. The state $x_i$ of a cell $i$ is - by hypothesis - governed by a continuous, deterministic and autonomous  dynamical system $\Phi_i$, on time $t \in \mathbb{R}, \ t \geq 0$, which we call  \em   the free dynamics of the cell $i$. \em It evolves on a finite-dimensional, differentiable, riemannian and compact manifold $X_i$.   Therefore, $x_i$ is a (maybe multidimensional) variable living on a compact   metric space with a local euclidean structure, such that for each instant $t \geq 0$ and for each initial state $x_i(0)$, the  continuous mapping $\Phi_i$ satisfies the following equalities:

%\begin{equation}
% ...
%\label{eqnprovisorio} \end{equation}

\begin{equation}  \begin{split}
\label{eqnPhi_i}  x_i(t) = \Phi_i(x_i(0), t), \ \ x_i(0) =\Phi_i (x_i(0),0), \ \  \\ x_i(t+s) = \Phi_i(x_i(0), t+s) = \Phi_i(\Phi_i(x_i(0), t), s) = \Phi_i(x_i(t), s) \ \forall \ t,s \geq 0.\end{split}\end{equation}

\noindent The flow defined on a finite dimensional manifold by an autonomous ordinary differential equation is an example of a continuous dynamical system, and may   model   each  cell. Each cell $i$ has its  own free dynamics, governed by its own rules. The spaces $X_i$ where the respective states $x_i$   of the cells live, and their dimensions,   may be mutually different.

We   call the scale of the whole network  the \em macro-scale \em or also, the \em global scale. \em It is defined by the  interactions $\Delta_{i,j}$ among  any ordered pair  $(i,j), \ i \neq j$ of different cells. These interactions    are produced, in our case, according to certain deterministic rules that we will state in Definition \ref{definitionSpikeMilestoneSatisfactionVariable}.

Each network ${\mathcal N}$ may be a cell of a larger hyper-macro network. We will not study this hyper-macro system, nor the role of ${\mathcal N}$ as a cell of it. Nevertheless, we assume that this hyper-macro system exists and may perturb   the state $x_i$ of   each individual cell $i$ in ${\mathcal N}$.

   Two main concepts are the \em spikes  or milestones \em and the eventual \em death \em of a cell.

\begin{definition} \em   
\label{definitionSpikeMilestoneSatisfactionVariable}
{\bf (Spikes or milestones, satisfaction  $S_i$ and goal  $\theta_i$)}

The \em spikes \em   or \em milestones \em of a cell $i$  are the instants $t_i$ when $i$ sends actions to the other cells $j $ of the network   ${\mathcal N}$. The spiking instants or milestones of the cell $i$   are defined by $i$, according to the value of a real variable $S_i$, which depends on the (maybe multidimensional) state $x_i$ of the cell. Thus, $S_i$ depends on time $t$, i.e. $S_i = S_i(x_i(t))$. We call $S_i$ \em the satisfaction variable \em of the cell $i$. This variable may be, for instance, the algebraic sum of several positive or negative components that depend  on the state $x_i$.

By hypothesis, the satisfaction variable $S_i$, while not perturbed from the exterior of the cell $i$, satisfies  the following two conditions,   (\ref{eqnSelfSatisfactionDiffEq}) and (\ref{eqnSelfSatisfactionResetValueX}):
 \begin{equation}
  \label{eqnSelfSatisfactionDiffEq}
 \frac{dS_i}{dt} = g_i (x_i), \ \ \ \mbox{ if } 0 \leq S_i \leq \theta_i,   \end{equation}  where \begin{equation}
  \label{eqnSelfSatisfactionDiffEq_bis} g_i \in C^1(X_i, \mathbb{R}^+), \ \ g_i(x_i) >0 \ \ \forall \ x_i \leq X_i   \ \ \ \mbox{ and } \ \ \  \theta_i >0. \end{equation}
  We call $\theta_i$  the \em goal level or threshold level\em.  It  may be achieved at many different states $x_i$ of the cell, i.e. the set $S_i^{-1}(\theta_i) \subset X_i$ is not necessarily reduced to a single point. Nevertheless, we assume, by hypothesis:
 \begin{equation} \label{eqnSelfSatisfactionResetValueX}
 \exists \ \mbox{ unique } x_{i, \mbox{\footnotesize reset } } \in X_i \ \  \mbox{ such that } \ \ S_i(x_{i, \mbox{\footnotesize reset } }) = 0.
 \end{equation}

 By hypothesis, at instant $0$ the initial state $x_i(0)$ of each cell satisfies:
 \begin{equation}\label{eqn87}0 \leq S_i(x_i(0)) < \theta_i.\end{equation}

 On the one hand,
 at each instant $t_i$ such that $\lim _{t \rightarrow t_i^-} S_i(x_i(t)) = \theta_{i}$, the cell $i$  reacts in the following way:

 First, the following \em reset or update rule \em holds:
 \begin{equation}
 \label{eqnSelfSatisfactionReset}
 \lim _{t \rightarrow t_i^-} S_i(x_i(t)) \geq \theta_{i} \ \ \Rightarrow \ \  x_i(t_i) = x_{i, \mbox{\footnotesize reset } } \ \Rightarrow \ \ S_i(x_i(t_i)) = 0, \end{equation} i.e. the state $x_i$ of the cell $i $ jumps to $x_{i, \mbox{\footnotesize reset } }$ at intant $t_i$. Thus, the satisfaction variable resets  to   zero at each milestone-instant.

  Second, the cell $i$ \lq\lq produces\rq\rq \  a  spike  at each milestone-instant $t_i$, i.e. $i$ plays in the game  ${\mathcal N}$ when its satisfaction variable $S_i$ arrives to the goal  $\theta_i$.  At instant $t_i$ the cell $i$ sends  instantaneous   signals $\Delta_{i,j}$   to   the other cells $j \neq i$.

On the other hand, when one and only one cell $j \neq i$ is spiking at an instant $t_j$ ($t_j$ may be different from $t_i$), it sends to $i$ a   signal $\Delta_{j,i} \geq 0$ (in general,
 $\Delta_{j,i} \neq \Delta_{i,j}$). In such a case, by hypothesis, the state $x_i$ suffers a discontinuity jump such that:
\begin{equation}\label{eqnInteraction0}\begin{split}S_i(x_i(t_j)) =   \ S_i (x_i(t_j^-)) +\Delta_{j,i} \ \ \mbox{ if } \ \ S_i (x_i(t_j^-)) +\Delta_{j,i} < \theta_i, \\    S_i(x_i(t_j)) = 0 \ \mbox{ otherwise.}\end{split}\end{equation}
(We denote $S_i(x_i(t_j^-)) = \lim_{t \rightarrow t^-_j} S_i(t)$.)

Finally, by hypothesis of the model, if all the cells of a nonempty set $I_n = \{j_1, \ldots, j_k\}$ - of $k$ different cells - spike simultaneously at some instant, say $t_n$, and if $\Delta_{j, i} \geq 0$ for all $j \in I_n$ and for all $i \not \in I_n$, then the state $x_i$ of any other cell $i \not \in I_n$ suffers a discontinuity jump such that:
\begin{equation}\label{eqnInteraction}\begin{split}
S_i(x_i(t_n)) =   \ S_i (x_i(t_n^-)) +\sum_{j \in I_n} \Delta_{j,i} \ \ \ \mbox{ if } \ \ S_i (x_i(t_n^-)) +\sum_{j \in I_n} \Delta_{j,i}< \theta_i, \\        S_i(x_i(t_n)) = 0 \ \mbox{ otherwise.}\end{split}\end{equation}

\end{definition}
\begin{definition} \em
\label{definitioncooperative} {\bf (Cooperative and antagonist cells)}

A cell $j$ is:

  \em   cooperative  \em if $\Delta_{j,i} \geq 0$ for all $i \neq j$ and $\max_{i} \Delta_{j,i} >0$,

   \em  fully cooperative  \em if $\Delta_{j,i} > 0$ for all $i \neq j$,

   \em  antagonist  \em if $\Delta_{j,i} \leq 0$   for all $i \neq j$ and $\min_i \Delta_{j,i} < 0$,

  \em fully antagonist  \em if  $\Delta_{j,i} <0$  for all $i \neq j$,

     \em mixed \em if   there exist $i_1 \neq j$ and $i_2 \neq j$ such that $\Delta_{i_1, j} >0$ and $\Delta_{i_2, j} < 0$.

\vspace{.3cm}

We call  a network \em  (fully) cooperative  \em if all its cells   are (fully) cooperative.

 \end{definition}

Formulae (\ref{eqnInteraction0}) and (\ref{eqnInteraction}) show  that  when each    cooperative cell $j$ spikes, then   it contributes to enlarge    the  values of the satisfaction variables $S_i$ of   the  other cells $i \neq j$. Thus, $j$ helps the other cells $i $ to approach to their respective goal levels, and so, it shortens the waiting times until the milestones of the others occur.

When an antagonist cell $j$ spikes, it reduces the values of the satisfaction variables $S_i$ of    the other cells Thus, $j$ enlarges the waiting times of the others to arrive to their respective goal levels.

In this paper, we will focus on full cooperative networks or subnetworks.

\vspace{.3cm}

By conditions (\ref{eqnSelfSatisfactionDiffEq}) and (\ref{eqnSelfSatisfactionDiffEq_bis}), while the free dynamics of a cell is not negatively perturbed by external agents, and while the satisfaction variable $S_i$ of the cell $i$ does not arrive to its goal level, it   strictly increases on time. In other words, each cell is \lq\lq born optimist\rq\rq: if free, it approximates with positive velocity to its goal. Besides,   the instantaneous velocity $g_i(x_i)$ is bounded away from zero, because $\min_{x_i \in X_i} g_i(x_i) >0$ exists, due to the compactness of the space $X_i$ and the continuity of $g_i$. So, if free, each cell arrives to its goal after a finite time. Nevertheless, if a negative term $- \delta_i < 0$ is added to the real function $g_i$, the velocity $dS_i/dt$  will decrease, and may also become negative.

\begin{definition} \em
\label{definitionInterference} {\bf (Negative external interferences)}

We call a negative number $- \delta  < 0$ \em negative differential interference  \em to the cell $i$ from its external environment, if during  an interval of time, the differential equation (\ref{eqnSelfSatisfactionDiffEq}) is substituted by
\begin{equation}
\label{eqnSelfSatisfactionDiffEq+delta}
\frac{dS_i}{dt} = g_i(x_i) - \delta  \ \ \mbox{ if }  0 \leq S_i < \theta_i.
\end{equation}

We call a negative number $-\Delta <0$ \em negative impulsive interference \em to the cell $i$ from its external environment, if at some instant $t$ for which the satisfaction variable $S_i$ has not arrived to its goal level, the state $x_i$ of the cell suffers a discontinuity jump such that
$$S_i(x_i(t)) =   S_i(x_i(t^-)) - \Delta.$$

\end{definition}

\begin{definition} \em
\label{definitiondeath} {\bf (Death of a cell)}

The \em death \em of a cell $i$ (cf. \cite{ErmentroutKopell1990}) occurs at an instant $T \geq 0$, if for any time $t >T$ the cell $i$ does not spike. In other words, after a cell $i$ dies, it does not arrive to its goal anymore, and so, it stops sending actions to the other  cells of the network forever.

In   Definition  \ref{definitionIntrinsicRisk}, we will state a  numerical formula to measure the theoretical intrinsic risk $R_i$ of death of any cell under eventual negative interferences, if it were not connected to the cooperative network.  In Definitions   \ref{definitionNetRisk} and \ref{definitionProtectionFactor}, we pose   formulae to measure the net risk of death $R'_i < R_i$, if  the cell $i$ is interacting in a cooperative network, and thus, the protection factor $P_i >0$ that the network provides to  $i$.
\end{definition}

\begin{definition} \em  \label{definitionGenericNetworks}

{\bf (Space of parameters)}

Let ${\mathcal N}$ be a network with $m \geq 2$ fixed cells. The \em  parameters of the network \em are:
\begin{equation}\label{eqn33} \mbox{Param}({\mathcal N}):= \left(\Big\{\big(\theta_i,  g_i\big)\Big\}_{1 \leq i \leq m}, \ \Big\{ \Delta_{i,j}  \Big \}_{1 \leq i , j \leq m, \ i \neq j} \right),\end{equation}
where, according to   Definition \ref{definitionSpikeMilestoneSatisfactionVariable},  $g_i $ is the given $C^1$ real function in the second member of the differential equation (\ref{eqnSelfSatisfactionDiffEq})   which governs the free dynamics of the cell $i$; $\theta_i \in {\mathbb{R}}^+$ is its   goal level; and $\Delta_{i,j} \in \mathbb{R}$ is the   interaction in the network  from the cell $i$ to the cell $j \neq i$.

In the space of all the   parameters of a network  with exactly $m$ cells (for  $m \geq 2$ fixed), we define the topology induced by the following metric (distance), where the parameters of the networks ${\mathcal N}$ and ${\mathcal N}'$ are defined by Equality (\ref{eqn33});   those denoted without ${}'$ (with ${}'$) correspond to the network ${\mathcal N}$ (resp. ${\mathcal N}'$):
\begin{equation}\label{eqnDistanciaParametros}\mbox{dist}\Big(\mbox{Param}({\mathcal N}), \ \mbox{Param}({\mathcal N}')\Big) := $$ $$ \max\Big \{\max_{1 \leq i \leq m}\big \{  |\theta_i - {\theta'}_i|, \ \|g_i - g'_i\|_{C_1}   \big \}, \  \max_{i \neq j}   |\Delta_{i,j} - \Delta'_{i,j}|    \ \Big\}.\end{equation}
 In the above equality,  $\| \cdot \|_{C^1}$ is the $C^1$ norm in the space of all the   $C^1$ real functions  defined on the compact manifold $X_i$.

\vspace{.3cm}

We say that a phenomenon of the global dynamics is \em robust, or persistent, or structurally stable \em if the set ${\mathcal G}$ of parameters of the networks for which the phenomenon occurs is open.
In other words, if a phenomenon is robust  and  if Param$({\mathcal N})  \in {\mathcal G}$, then Param$({\mathcal N'})$ still belongs to    ${\mathcal G}$ for any sufficiently small perturbation  ${\mathcal N}'$ of the network ${\mathcal N}$ (openness  condition).

\end{definition}

\begin{definition} \em  {\bf (Spiking instants and interspike intervals of the network)}

We denote by $\{t_n\}_{n \in \mathbb{N}^+}$ the strictly sequence of all the instants for which at least one neuron spikes. We call $t_n$ the \em $n$-th. spiking instant \em of the network. We call $(t_n, t_{n+1}]$ the \em $n$-th. interspike interval \em of the network and denote it by $ISI^{(n)}$.

\end{definition}
\begin{definition} \em {\bf (Clusters or spiking codes)} \label{definitionCluster}

We denote by $I_n$ the set of neurons that spike at the instant $t_n$. We call $I_n$ the \em $n$-th. cluster  \em or also, the \em $n$-th. spiking code. \em
\end{definition}

\vspace{.5cm}

\section{Mathematical statements} \label{sectionStatementOfResults}

\begin{definition} \em  {\bf (Spike-synchronization)}
 \label{definitionSincronizacion}

We say that the network \em eventually synchronizes spikes \em   if there exists a subsequence $\{t_{n_h}\}_{h \in \mathbb{N}}$ of  spiking instants    such that the respective clusters $I_{n_h} $ are $ \{1, \ldots, m\}$ (i.e. all the cells spike at instants $t_{n_h}$).

\end{definition}

\begin{definition} \em
{\bf (Periodic spike-synchronization)} \label{definitionSincronizacionPeriodica}

We say that all the cells of the network \em eventually periodically synchronize  spikes with period $p \geq 1$, \em   if  there exists $n_0 \geq 0$ such that:

\vspace{.2cm}

\noindent {\bf i) } the subsequence $\{t_{n_h}\}_{ h \in \mathbb{N}}$ of Definition \ref{definitionSincronizacion} satisfies
$$t_{n_h} = t_{n_0 + hp} \ \ \forall \ h \in \mathbb{N},$$

\noindent {\bf ii) } the sequence  $\{I_n\}_{n \geq 0}$ of clusters satisfies
$$I_n = I_{n+ p} \ \ \forall \ n \geq n_0,$$

\noindent {\bf iii) }  the sequence  $\{  t_{n+1} - t_n\}_{n \geq 0} $ of the   interspike intervals' lengths satisfies:
$$t_{n+ p + 1} - t_{n+p} = t_{n+1} - t_{n} \ \ \forall \ n \geq n_0.$$

\vspace{.3cm}

 We call $p$ the \em natural spiking period  \em of the network.
\end{definition}
Note that    the network eventually periodically synchronizes spikes, if and only if at least one    cell  spikes at instant $t_n$ for all $n \geq 0$,   no cell spikes at instants $t \in (t_{n}, t_{n+1})$ for all $n \geq 0$, and all the cells spike at instants $t_{n_0 + hp}$ for any natural number  $h \geq 0$.

\begin{definition} \em  {\bf (Transitory time)}

 For each fixed initial state for which the network periodically synchronize spikes, we call $t_{n_0}$ the \em waiting time \em or \em transitory time \em until the synchronization of the full network occurs.

 Note that the occurrence of synchronization and the value of the waiting time  depend  on the initial state of the network.

\end{definition}

Let ${\mathcal N}$ be a network   composed    by $m \geq 2$  fully  cooperative cells   (cf. Definition \ref{definitioncooperative}). Recall that $\theta_j$ denotes the    goal level of the cell $j$, and $\Delta_{i,j}$ denotes the action from the cell $i$ to the cell $j$  (cf. Definition \ref{definitionSpikeMilestoneSatisfactionVariable}).

\begin{definition} \em
\label{definitionLargeCooperativeNetwork} {\bf (Large cooperative network)}

A fully cooperative network is called \em large \em if
\begin{equation}
\label{eqnm'>Theorem1}
\sqrt {m \ } > \max\Big\{\sqrt 3, \ \ \  \frac{\displaystyle \max \{\theta_j \colon   \ j \in {\mathcal N} \}}{\displaystyle \min\{ \Delta_{i,j}\colon  \  i, j \in {\mathcal N}, \ i \neq j\}} \ + \ 1\Big\},\end{equation}
\end{definition}

%\textcolor{red}{quitar newpage}

%\newpage

 \label{subsectionStatementOfResults}

\begin{theorem} {\bf (Synchronization) }  \label{theoremFullCooperativism}

If ${\mathcal N}$ is a large fully cooperative network
   then:

\noindent{\bf (a) } From any initial state, all the cells  of the network ${\mathcal N}$  eventually periodically synchronize spikes.

\noindent{\bf (b) } The eventual periodic synchronization  is a robust phenomenon.

\end{theorem}

We prove Theorem \ref{theoremFullCooperativism} in Subsection  \ref{proofPart(a)Theorem1}.

\begin{theorem} {\bf (Upper bounds for the transitory time and  the spiking period)}
\label{TheoremPart(c)}

Under the hypothesis of Theorem \em \ref{theoremFullCooperativism}, \em the transitory time $T $ and the natural spiking period $p$  satisfy the following inequalities:

\begin{equation} \label{eqnBoundTransitoryTime} T \leq    \max_{1 \leq i \leq m}   \Big \{ \frac{  \ \ \ \ \theta_i   \ \ \ \ }{\displaystyle  \ \min \ \{g_i(x_i) \colon \ x_i \in X_i \   \} } \ \Big \},\end{equation}
\begin{equation} \label{eqnBoundPeriodp}p \leq  1 + \frac{\displaystyle \max \{\theta_j \colon   \ j \in {\mathcal N} \}}{\displaystyle \min\{ \Delta_{i,j}\colon  \  i, j \in {\mathcal N}, \ i \neq j\}}.\end{equation}

\end{theorem}

We prove Theorem \ref{TheoremPart(c)} in Subsection \ref{proofPart(c)Theorem1}.

\vspace{.3cm}

Now, we   define   the risk of death of each cell and the protection factor of the network.

In Definition \ref{definitiondeath} we say that a cell $i$  dies
if (due to external causes of the cell) there is an instant  $T$
such that for all $t > T$ the satisfaction variable $S_i$ remains
smaller than the goal level $\theta_i$. In brief:
$$ i \mbox{ dies } \ \Leftrightarrow \ \ S_i(x_i(t)) < \theta_i \ \ \forall \ t > T  \mbox{ for some } T \geq 0.$$ Due to the rules of the interactions among the cells, if a cell $i$ dies, it can not
act on the network anymore (after time $T$).

We measure the     risk of death     of the cell $i$ by the following definitions:
%\subsection{Previous definitions and lemma}

\begin{definition} \em
\label{definitionIntrinsicRisk}
{\bf (Intrinsic risk of death $R_i$)}

The \em intrinsic risk of death $R_i$ \em  of the cell $i$, relatively to the other cells of the network is defined by

\begin{equation}\label{eqnInstrinsicRisk}R_i := \frac{ \theta_i}{\displaystyle \max_{1 \leq j \leq m} \theta_j} ,\end{equation}
where $m$ is the number of cells of the network.
 Since $\theta_i >0$ for all $i$, Equality (\ref{eqnInstrinsicRisk}) immediately implies the following:
 \begin{equation}\label{eqnR_i>0}0  < R_i \leq 1. \end{equation}
\end{definition}

Recall Definition \ref{definitionSpikingCodes}: the $n$-th. spiking-code   $I_n$ is the (nonempty) set of cells that spike at instant $t_n$. Let us fix a cell $i$.

\begin{definition} \em {\bf (Spiking instants of each cell)}.

The  sequence $\{t_i^{(h)}\}_{h \geq 1}$ of \em spiking times \em of the cell $i$ is defined by the following equalities:
$$\begin{array}{lcll}
    t_i^{(1)} &:= &\min\{t_n >0: \ \ i \in I_n\} &\mbox{ if } i \in I_n \mbox{ for some } n \geq 0, \\  t_i^{(1)} &:= &+ \infty & \mbox{ otherwise},  \\
t_i^{(h+1)} &:= &\min  \{  t_n > t_i^{(h)}: \ \ i \in I_n\} &\mbox{ if } i \in I_n \mbox{ for some } t_n >  t_i^{(h)}, \\ t_i^{(h+1)} &:=& + \infty &\mbox{ otherwise}.
\end{array}$$
We denote
$$t_i^{(0)} := 0.$$

For each cell $i$ and for each natural number $h \geq 1$ such that $t_i^{(h)} < + \infty$, we call $t_i^{(h)}$ \em the $h$-th. spiking instant or the $h$-th. milestone of the cell $i$. \em
Note that  a cell $i$ dies   if and only if $t_i^{(h+1)} = + \infty$ for some $h \geq 0.$
\end{definition}
\begin{definition} \em
\label{definitionISI_i} {\bf (The inter-spike-intervals  of each cell)}.

We call the time-interval $$\mbox{ISI}^{(h)}_i :=     {  \big (} t_i^{(h)},t_i^{(h+1)}{  \big]} \ \ \forall \ h \geq 0 $$ \em the $h$-th. inter-spike-interval of the cell $i$. \em Note that we include the instant $t_i^{(h+1)}$  in the $h$-th. inter-spike-interval of the cell $i$.
\end{definition}

When a cell $i$ receives  cooperative actions $\Delta_{j,i}>0$ from the other cells $j \neq i$ of the network, then its satisfaction variable $S_i$ is increased, to approach (or even reach) its goal $\theta_i$. Equivalently, a positive action $\Delta_{j,i} >0$ can be understood, as a reduction of the goal level $\theta_i$ substituting it by $\theta_i - \Delta_{i,j}$.  Roughly speaking, when a cell $i$ receives positive interactions from the other cells, its risk of death   diminishes.
\begin{definition} \em
For each fixed cell $i$, and for each fixed natural number $h \geq 0$, we denote:
\begin{equation}\label{eqnNetDelta}\Delta_{i, {\mathcal N}}^{(h)}:= \sum_{\displaystyle j \neq i, \ j \in I_n, \ t_n \in \mbox{\footnotesize ISI}_i^{(h)}} \Delta_{j,i}, \end{equation}

We call $\Delta_{i, {\mathcal N}}^{(h)}$ \em the   net  action \em that the cell $i$ receives  during its $h$-th. inter-spike-interval  from the other cells of the network.
\end{definition}

If the network is not   cooperative, then some of the actions $\Delta_{j,i}$ in   Formula (\ref{eqnNetDelta}) may be negative, so the sum of all of them may be negative or null. Note that, in general, the net action $\Delta_{i, {\mathcal N}}^{(h)}$ depends on the initial condition of the network.

\begin{definition} \em
\label{definitionNetRisk} {\bf (Net risk of death $R'_i$)}

The \em net risk of death ${R'}^{(h)}_{i, {\mathcal N}}$ \em of the cell $i$ while connected to the network ${\mathcal N}$, during its $h$-th. inter-spike-interval, is
\begin{equation} \label{eqnNetRisk} {R'}^{(h)}_{i, {\mathcal N}} := \max\Big \{0, \  \min\Big\{1, \   \frac{\displaystyle  \ \ \theta_i - \Delta_{i, {\mathcal N}}^{(h)} \ \ }{\displaystyle \max_{1 \leq j \leq m} \theta_j} \Big\} \Big \}.\end{equation}
\end{definition}

In general the net risk of death ${R'}^{(h)}_{i, \mathcal N}$ depends on the initial condition. For simplicity in the notation, in the sequel we will write $R'_i$ to denote the net risk of the cell $i$, when the network ${\mathcal N}$ is clear from the context, and when referring to some fixed $h \geq 0$.

\begin{definition} \em
\label{definitionProtectionFactor} {\bf (Protection factor $P_i$)}
The \em protection factor $P_{i, {\mathcal N}}^{(h)}$ \em of the network ${\mathcal N}$ to its cell $i$, during the $h$-th. inter-spike-interval of   $i$, is
 \begin{equation} \label{eqnProtectionFactor} P_{i, {\mathcal N}}^{(h)} :=  \frac{ \ \  \displaystyle   \Delta_{i, {\mathcal N}}^{(h)} \ \ }{\theta_i} .\end{equation}
\end{definition}
In the sequel, we will simply write $P_i$ instead of $P_{i, {\mathcal N}}^{(h)}$, when the network ${\mathcal N}$ is clear from the context, and when referring to some fixed $h \geq 0$.

\begin{remark} \em 
\label{remarkProtectionFactor}

From Equality (\ref{eqnProtectionFactor}), it is immediate to deduce that the protection factor $P_i$ is positive (negative) if and only if the minimum net sum of actions $\Delta_{i, {\mathcal N}}^{(h)}$ that the cell $i$ receives from the other cells of the network is positive (resp. negative). Therefore, if the network is fully cooperative,   its protection factor is   strictly positive, and if the network is fully antagonist, its protection factor is strictly negative.

\end{remark}

In the following Proposition \ref{propositionProtectionFactor} we state that net risk $R'_i$ of a cell that is connected to the network ${\mathcal N} $, is (essentially) the product of its intrinsic risk $R_i$, if $i$ were not connected to the network, by   $1- P_i$.

In other words, if the protection factor $P_i$ of the network were positive, then the net risk $R'_i$ is smaller than the intrinsic risk of death, and if the protection factor $P_i$ were 1, then the net risk of death $R'_i $ would be zero. So, a cooperative network always reduces the intrinsic risk, and may also make it zero, if the cooperative interactions are large and frequent enough.

\begin{proposition}
\label{propositionProtectionFactor} {\bf(Formula of the protection factor of the network)}

 \noindent  \em
\begin{equation} \label{eqnFormulaProtectionFactor} R'_i := \max \Big \{0, \ \min\big\{1, \ \ \big(1 - P_i \big) \ R_i \big \} \Big\}.\end{equation}
\end{proposition}
We prove Proposition \ref{propositionProtectionFactor} in Subsection \ref{subsectionProofPrpositionProtectionFactor}.

\begin{theorem} {\bf (Protection factor)}
\label{TheoremPart(d)}

\noindent  Under the hypothesis of Theorem \em \ref{theoremFullCooperativism}, \em the protection factor $P_i$ of the network  to each cell $i$  is positive. So, the net risk of death $R'_i$  of each cell $i$ under negative external interferences, when it is connected to the network,  is strictly smaller than its intrinsic risk $R_i$ when it is isolated from the network.

\end{theorem}

We prove Theorem \ref{TheoremPart(d)} in Subsection \ref{proofPart(d)Theorem1}.

\begin{definition} \em
{\bf (Code-patterns) } \label{definitionSpikingCodes}

For any natural number $k \geq 1$, and for any fixed                             initial state $${\mathbf x}(0) = (x_1(0), \ldots, x_m(0))$$                      of the network, we construct the following word of clusters                      with length $k$:
$$\Pi_{n,k} := (I_n, \ldots, I_{n+k-1}).$$
We call $\Pi_{n,k}$ the \em $n-$th. code-pattern with                            length $k$ from the initial state ${\mathbf x}(0)$. \em

\end{definition}
\begin{definition} \em  {\bf (Recurrent code-patterns)}

We say that a   code-pattern $\Pi_{k}$ with length $k$  is \em recurrent, \em if there exists $n_j \rightarrow + \infty$ such that
$$\Pi_{n_j,k} = \Pi_k \ \ \forall \ j \in \mathbb{N}.$$
We denote by ${\mathcal P}_k$   the set of all the recurrent code-patterns with length $k$, obtained from all the initial states of the network. We denote by
$$\#  {\mathcal P}_k \geq 1$$
  the cardinality of the set ${\mathcal P}_k$.
\end{definition}

\begin{definition} \em
\label{definitionEntropy}  {\bf (Amount of information)}

We denote
\begin{equation} \label{eqnDefinitionEntropia}H   := \sup_{k \geq 1} \{\log_2 (\#{\mathcal P}_k) \} \ \mbox{bits}  \in [0, + \infty],\end{equation}
where $\log_2$ is the   logarithm in base 2.

We call $H $ \em the total  amount of information \em that the network ${\mathcal N}$ is able to process.
\end{definition}

Note that, since $\#{\mathcal P}_k  \geq 1$ for all $k $, then $H  \geq 0$.

\begin{definition} \em
{\bf (Entropy)}

If  $H= + \infty$, we define \em the entropy   $h$, \em or the exponential rate of increasing information, by
$$h := \limsup_{k \rightarrow + \infty} \frac{H_k}{k}, \ \mbox{ where } \ H_k := \log_2 (\#{\mathcal P}_k).$$

\end{definition}

\noindent {\bf Interpretation of the amount of information $H$. } The  number $\#{\mathcal P}_k$ of different recurrent patterns $\Pi_k$ with length $k$   measures the dynamical richness of ${\mathcal N}$, with respect to the many possible finite words that the spiking-codes can show. Namely, $H$ takes into account how many different   clusters $I_n$  are exhibited  at the milestones $t_n$ of the network.   Thus, the total     amount of information $H $ is a quantitative index of the maximum  global codified dynamical diversity  that the network ${\mathcal N}$    is able to show along   its recurrent evolution in the future,  if the observer   looks at  \em   each cell with an individual role. \em In fact, note that  if $i \neq j$, then the event for which the cell $i$ spikes and $j$ does not, is \em distinguished \em    from the converse event.

\begin{theorem} {\bf (Amount of information)}
\label{TheoremPart(e)}

 Under the hypothesis of Theorem \em \ref{theoremFullCooperativism}, \em
the total amount of information $H$ of the network   is
$$H = \log_2 p \leq \log_2 \Big(1 + \frac{\displaystyle \max \{\theta_j \colon   \ j \in {\mathcal N} \}}{\displaystyle \min\{ \Delta_{i,j}\colon  \  i, j \in {\mathcal N}, \ i \neq j\}} \Big ).$$

\end{theorem}
We prove Theorem \ref{TheoremPart(e)}    in Section \ref{proofPart(e)Theorem1}.
\begin{definition} \em
{\bf (Mutually similar cells)}

We say that   the cells of a fully cooperative network are \em mutually similar \em (with respect to the minimum interaction) if:

\begin{equation}
\label{eqnhypothesisnueva}
\frac {\displaystyle\Big(\min_{1 \leq i \leq m}   \theta_i \Big) \Big(\min_{1 \leq i \leq m}\min\{g_i (x_i) \colon x_i \in X_i\}\Big)} {\displaystyle\Big( \max_{1 \leq i \leq m}   \theta_i \Big) \Big(\max_{1 \leq i \leq m}\max\{g_i (x_i) \colon x_i \in X_i\}\Big) }> 1 - \frac{  \displaystyle \min \{  \Delta_{i,j} : \ i,j \in {\mathcal N}, \ i \neq j\}  }{\displaystyle\max_{1 \leq i \leq m} \theta_i },
\end{equation}
\end{definition}

\begin{theorem} {\bf (Cooperative networks of mutually similar cells)}
\label{TheoremParts(fgh)}
 Under the hypothesis of Theorem \em \ref{theoremFullCooperativism}, \em
if besides  the cells are mutually similar, then

 \vspace{.2cm}

 \noindent {\bf (a) } From any initial state all the cells eventually periodically synchronize all their spikes  with natural spiking period $p= 1$.

\vspace{.2cm}

\noindent{\bf (b) } After the transitory time $T$ has elapsed, the protection factor $P_i$ of the network  to each cell $i$  is 1. So, its net risk of death $R'_i$  under negative external interferences  is null.

\vspace{.2cm}

\noindent{\bf (c) }
The total amount of information $H$ of the network  is zero
\end{theorem}

\vspace{.2cm}

We prove Theorem \ref{TheoremParts(fgh)} in Section \ref{proofPart(f)Theorem1}.

\vspace{.2cm}

\noindent{\bf Remark: } In particular Inequality (\ref{eqnhypothesisnueva}) - and  thus, also the assertions (a), (b) and (c) of Theorem \ref{TheoremParts(fgh)} -  holds  if all the cells are mutually identical, i.e. $$\theta_i = \theta_j = \theta, \ \ g_i = g_j = g: X \mapsto \mathbb{R}^+  \ \ \forall \  i \neq j,$$  and if   the minimum positive interaction is large enough, i.e.:
 $$\frac{  \displaystyle \min \{  \Delta_{i,j} : \ i,j \in {\mathcal N}, \ i \neq j\}  }{\theta } > 1 -
 \frac {   \min\{g  (x ) \colon x  \in X \} } { \ \max\{g  (x ) \colon x  \in X \} }.$$

\vspace{.2cm}

Now let us pose a result about large cooperative networks whose graphs are not necessarily complete.

\begin{definition} \em
{\bf (Full cooperative core)} \label{definitionCore}

Let ${\mathcal N}$ be a network. A \em full cooperative core \em in ${\mathcal N}$, if it exists, is a subnetwork   ${\mathcal N}_1$ composed by   fully cooperative cells, i.e. $$i \in {\mathcal N}_1 \ \Rightarrow \ \ \Delta_{i,j} > 0 \ \forall \ j \in {\mathcal N}, $$ and such that  all the cells that do not belong to ${\mathcal N}_1$ have non negative actions on the cells of ${\mathcal N}$, i.e.
$$i \not \in {\mathcal N}_1 \ \Rightarrow \ \Delta_{i,j} \geq 0 \ \forall \ j \in {\mathcal N}.$$

\end{definition}
\begin{corollary}
\label{corollary1}   {\bf (Networks with a large     cooperative core.)}

  Let ${\mathcal N}$ be a network that has a full cooperative core ${\mathcal N}_1$. Assume that the number $m$ of cells in ${\mathcal N}_1$ is large enough. Precisely:
  \begin{equation} \label{eqn8cor1} \sqrt {m \ } > \max\Big\{\sqrt 3, \ \ \  \frac{\displaystyle \max \{\theta_j\colon \displaystyle j \in {\mathcal N}\}}{\displaystyle \min \{ \Delta_{i,j}\colon \displaystyle i \neq j, \ i \in {\mathcal N}_1, \ j \in {\mathcal N}\}} \ + \ 1 \ \Big\}.\end{equation}
  Then Assertions   of Theorems \em \ref{theoremFullCooperativism}, \ref{TheoremPart(c)}, \ref{TheoremPart(d)}, \ref{TheoremPart(e)} \em and \em \ref{TheoremParts(fgh)} \em hold, where \em
  $\displaystyle \min\{ \Delta_{i,j}\colon  \  i, j \in {\mathcal N}, \ i \neq j\} $ \em must be substituted by \em $ \displaystyle \min\{ \Delta_{i,j}\colon  \  i \in {\mathcal N}_1, \  j \in {\mathcal N}, \ i \neq j\}.$
\end{corollary}

We prove Corollary \ref{corollary1} in Section \ref{sectionProofOfCorollaries}.

\vspace{1cm}

\section{The proof of synchronization} \label{sectionSincronizacion}

In this section we   prove   Theorems \ref{theoremFullCooperativism} and \ref{TheoremPart(c)}.
To do so, we need the following previous result:

\begin{lemma}
\label{propositionSynchronizationPeriodic}

\vspace{.2cm}

\noindent {\bf (A) } If from some initial state of the network  there are  two different (minimal) instants $0 < t^*_0 < t^*_p$ when all the cells spike simultaneously, then the network eventually periodically synchronizes spikes with some natural spiking period $p \geq 1$.

\noindent {\bf (B) } If from all the initial states the network eventually periodically synchronizes spikes with the same period $p \geq 1$, then  the number of all the possible recurrent code-patterns with length $k \geq 1$ is
$$\#{\mathcal P}_k\leq p \ \ \ \forall \ k \geq 1,$$
and
$$\max_{k \geq 1} (\#{\mathcal P}_k) = p.$$

\noindent {\bf (C) } If from some initial state the network eventually periodically synchronizes  spikes with period $p= 1$, then,  the sequence $\{I_n\}_{n \geq n_0}$ of   clusters is constantly equal to the set of all the cells of the network, i.e.  $$I_n = \{1, 2, \ldots, m\} \ \ \forall \ n \geq n_0.$$

\noindent {\bf (D) } If from all the initial states the network eventually periodically synchronizes   spikes with period $p= 1$, then  the number of all the possible recurrent code-patterns  with length $k \geq 1$ is
$$\#{\mathcal P}_k= 1 \ \ \ \forall \ k \geq 1.$$
\end{lemma}

\noindent{\em Proof: }
Assertions (C) and (D) are immediate consequences of (A) and (B) respectively,   in the particular case $p= 1$.

Let us prove Assertion (A). Recall that $m \geq 2$ denotes the number of cells in the network.

By hypothesis, all the cells of the network spike at instants $t^*_{0} > 0 $  and $t^*_p > t^*_0$.   Thus, due to the reset hypothesis in Definition \ref{definitionSpikeMilestoneSatisfactionVariable}, Formula (\ref{eqnSelfSatisfactionReset}), the state of the network
at instant $t^*_{0}$,    is $${\mathbf x}(t^*_0) = {\mathbf x}_{\mbox{\footnotesize reset}} := \big(x _{1, \mbox{\footnotesize reset}}, \ldots, x _{i, \mbox{\footnotesize reset}}, \ldots, x _{m, \mbox{\footnotesize reset}}    \big ),  $$ where $x _{i, \mbox{\footnotesize reset}}$ is the unique   point in the state-space $X_i$ that satisfies  Equality (\ref{eqnSelfSatisfactionResetValueX}). Thus $${\mathbf S}({\mathbf x} (t^*_{0})) = {\mathbf 0}, \ \ \mbox{ i.e. }
 S_i(x_i(t^*_{0})) = 0 \ \ \forall \ 1 \leq i \leq m.$$
As the dynamics defined by (\ref{eqnPhi_i}) is deterministic, once fixed the unique state ${\mathbf x}_{\mbox{\footnotesize reset}} $ of the network, a unique orbit exists in the future. Thus, translating the origin $0$ of time to $t^*_0$, and recalling Equalities (\ref{eqnPhi_i}), we obtain:
$${\mathbf S}({\mathbf x} (t)) = {\mathbf S}({\mathbf x} (t^*_0 + (t-t^*_0))= {\mathbf S}({\mathbf {\Phi}}({\mathbf {\Phi}} ({\mathbf x_0},t^*_0), t- t^*_0)) \ \ \forall \ t \geq t^*_0,$$
where ${\mathbf {\Phi}} = (\Phi_1, \ldots, \Phi_i, \ldots, \Phi_m)$.
Since $S_i(x_i(t))$ satisfies the     differential equation (\ref{eqnSelfSatisfactionDiffEq}) for all $t >t^*_0  $ such that $t \neq t_n$, and $${\mathbf S}( {\mathbf {\Phi}} ({\mathbf x_0},t^*_0)) = {\mathbf 0},  $$ we have:
$${\mathbf x}(t^*_0) ={\mathbf \Phi} ({\mathbf x}_0, t^*_0) = {\mathbf x}_{\mbox{\footnotesize reset}},$$
$$  {\mathbf S}({\mathbf x} (t)) = {\mathbf S}({\mathbf {\Phi}}({\mathbf {\Phi}} ({\mathbf x_0},t^*_0), t- t^*_0)) = {\mathbf S}  ({\mathbf {\Phi}} ({\mathbf x}_{\mbox{\footnotesize reset}}, t- t^*_0)) =  {\mathbf S}  ({\mathbf {x}}^* ( t- t^*_0)) \ \ \forall \ t \geq t^*_0,$$
where ${\mathbf x}^* = (x^*_1, \ldots, x^*_i, \ldots, x^*_m) $ is the unique solution  ${\mathbf x}^*(\cdot) = {\mathbf {\Phi}} ({\mathbf x}_{\mbox{\footnotesize reset}}, \cdot)$ of the deterministic dynamical system    (\ref{eqnPhi_i})  with initial state  ${\mathbf x}^*(0) =  {\mathbf x}_{\mbox{\footnotesize reset}}$, plus the deterministic interactions' rule (\ref{eqnInteraction})during the time-interval $(t^*_0, t^*_p)$.

By the hypothesis, all the cells spike simultaneously again at the instant  $t^*_p > t^*_0$. Therefore,
 $  {\mathbf S}({\mathbf x} (t^*_p)) = {\mathbf S}({\mathbf x} (t^*_0)) = {\bf 0}$. Thus  ${\mathbf x} (t^*_p)= {\mathbf x} (t^*_0) = {\mathbf x}_{\mbox{\footnotesize reset}}  $. Then,
 $${\mathbf x} (t^*_p) = {\mathbf x}^* ( t^*_p-t^*_0) =  {\mathbf x}_{\mbox{\footnotesize reset}} = {\mathbf x}^* ( 0).  $$ We deduce that  the unique solution ${\mathbf x}^*$  which has initial condition ${\mathbf x}_{\mbox{\footnotesize reset}}$, is periodic with period $t^*_p- t^*_0$. Thus, the instants $t_n $ and the spiking-codes $I_n$ are determined recursively from   the unique \em periodic \em orbit ${\mathbf x}^*$. We deduce  \begin{equation} \label{eqnPropositionSycnh1}I_n = I_{n+p} \ \ \forall \ n \geq n_0.\end{equation}
Also, \begin{equation} \label{eqnPropositionSycnh2}I_{n_0 + hp}= \{1, 2 \ldots, m\}, \ \   I_n \subset_{\displaystyle \neq } \{1, 2 \ldots, m\} \ \ \forall \ n_0 + hp <  n < n_0 + (h+1)p, \ \ \forall \ h \geq 0,
 \end{equation}
and the sequence $\{t_n\}_{n \geq 0}$ of instants for which at least one cell spikes   satisfies:  \begin{equation} \label{eqnPropositionSycnh3}
   t_{n+ p+1} - t_{n + p} = t_{n+1} - t_n \ \ \ \forall \ n \geq n_0.\end{equation}  Equalities (\ref{eqnPropositionSycnh1}), (\ref{eqnPropositionSycnh2}) and (\ref{eqnPropositionSycnh3}) end the proof of Assertion (A) of Lemma \ref{propositionSynchronizationPeriodic}.

 \vspace{.2cm}

 Let us prove Assertion (B):

  \noindent By Assertion (A), the sequence $\{I_n\}_{n \geq n_0}$ is periodic with period $p$. Thus, from Definition \ref{definitionSpikingCodes}, for any fixed natural number $k \geq 1$, all the recurrent code-patterns with length $k$ are:
 \begin{equation}\label{eqn03} \begin{split}
 (I_{n_0}, I_{n_0 +1}, \ldots, I_{n_0 + k-1}), \ (I_{n_0 +1}, I_{n_0 +2}, \ldots, I_{n_0 + k}), \ldots,  \\
 \ldots,   (I_{n_0 + r}, I_{n_0 +r+1}, \ldots,  I_{n_0 + r+k-1}), \ldots, \\   \ldots, (I_{n_0 + p-1}, I_{n_0 + p}, \ldots, I_{n_0 + p+k-2})\end{split}\end{equation}
 with $0 \leq r \leq p-1$. In fact, for $n= n_0 + hp+ r \geq p$ the code-pattern \\ $(I_n, I_{n+1}, \ldots, I_{n+k-1})$ coincides with $(I_{n_0 +r}, I_{n_0 +r+1}, \ldots, I_{n_0 +r+k-1})$, because $I_n = I_{n_0 + hp+ r} = I_{n_0 +r}$.
 So, all the code-patterns in the list (\ref{eqn03}) are recurrent. Two or more code-patterns in the list (\ref{eqn03}) may coincide, so the number of different code-patterns with length $k$ is at most equal to the number of items in the list
 (\ref{eqn03}). Thus:
 $$\#{\mathcal P}_k \leq p \ \ \forall \ k \geq 1.$$
 Now let us prove that, in the particular case that $k= p$, the code-patterns of the list (\ref{eqn03}) are pairwise different. To prove this assertion, with no loss of generality, we assume $n_0 = 0$ (if not, we translate the origin 0 of time to the instant $t_{n_0}$). Assume that
 \begin{equation}
 \label{eqn04}
 (I_r, I_{r+1}, \ldots, I_{r+p-1}) = (I_s, I_{s+1}, \ldots, I_{s+p-1})\end{equation}
 for $0 \leq r , s \leq p-1$. We must prove that $r= s$.

 The code-pattern $I_0= \{1, \ldots, m\}$ will appear once in both
   patterns (\ref{eqn04}), because they both have length $p$, which  is the (minimum) period
  of the sequence $\{I_n\}_{n \geq 1}$.
  Say $I_0 = I_{r+h} = I_{s + k}$  with $0 \leq h,k \leq p-1$.

  Both patterns in Equality (\ref{eqn04}) coincide. Then the positions  $h$ and $k$   are the same: $$h = k.$$ Besides,
  since $0 \leq r+h, s+k \leq p-1$ and $0 \leq r,s \leq p-1$,
  we have $$|(r+h)-(s+k)| = |r-s| \leq p-1< p.$$

  As $I_0 = I_{r+h} = I_{s + k}$, there are two indexes
  $r+h$ and $s+k$, whose difference is smaller than $p$,
   such that the   respective patterns coincide with
  $I_0 = \{1, \ldots, m\}$. In other words, all the cells spike at two instants $t_{r+h}$
  and $t_{s+k}$ such that  $|(r+h)-(s+k)|<p$. But $p$ is the minimum positive natural number such that all
  cells spike at instants $t_n$ and $t_{n+p}$  for some $n$. We deduce that $r+h = s+k$. Since
  we already know that $h = k$,   we deduce that $r= s$, as wanted.

  We conclude that all the code-patterns of the list (\ref{eqn03}) are pairwise different when
  the length $k$ equals the period $p$. Thus, the number
  of different code-patterns with length $p$ is $p$, ending the proof of Assertion (B) of Lemma \ref{propositionSynchronizationPeriodic}.
\hfill $\Box$

\subsection
{Proof of   Theorem \ref{theoremFullCooperativism}}
\label{proofPart(a)Theorem1}

\noindent{\em Proof: }
{\bf Part a) }

From Lemma \ref{propositionSynchronizationPeriodic}, to prove that all the cells of the network eventually periodically synchronize spikes, it is enough to prove that there exist two instants $0 < t_0 < t_p$ such that all the cells simultaneously spike at $t_0$ and at $t_p$.

If \em for any initial condition \em we find a single instant $t_0 > 0$ at which all the cells simultaneously spike, then, taking as new initial state $ {\mathbf x}(t_0)) = {\mathbf x}_{\mbox{\footnotesize reset}}$, we deduce that there exists a second instant $t_p > t_0$ at which all the cells simultaneously spike. Thus, to prove Part (a) of Theorem \ref{theoremFullCooperativism}, it is enough to show the following:

\vspace{.2cm}

 \noindent{ Assertion {\bf (i) } to be proved: }  \em For any initial condition, there exists an instant $t_0 >0$ at which all the cells simultaneously spike. \em

 \vspace{.2cm}

 From any   initial state ${\mathbf x}_0 = \big( x_1(0), \ldots, x_m(0) \big)$ such that $0 \leq S_i(x_i(0)) < \theta_i$ for all $1 \leq i \leq m$, consider the  state   $$ {\mathbf x }(t)  = \Big ( x_1(t), \ldots, x_i(t), \ldots,  x_m(t) \Big), $$ and  the $m$-dimensional vector whose components are   the satisfaction variables $S_i(x_i(t))$   at instant $t > 0$.

    Since each variable $S_i$ is governed by the differential equation (\ref{eqnSelfSatisfactionDiffEq}) with a strictly positive real function $g_i:X_i \mapsto \mathbb{R}^+$ (which is continuous on the compact space $X_i$), plus  the eventual sum of interactions $\Delta_{j,i} \geq 0$ according to Equalities (\ref{eqnInteraction}),   we deduce:

   \vspace{.2cm}

   \noindent{Property {\bf (ii) }}    \em While no interference from outside the network appears,  for each cell $i $  the real variable $S_i $  { is strictly increasing on} $t$,  for all $t \geq 0 $  {such that} $S_i(x_i(t)) \in \big [0, \theta_i  \big ) $. Besides, except at those instants where $S_i$ is discontinuous,  its derivative respect to time $t$ exists,  is   positive and   bounded away from zero, and at the instants where $S_i$ is discontinuous, the discontinuity jumps are positive.\em

\vspace{.2cm}

    We are assuming that $S_i(x_i(0)) < \theta_i$.
    Thus, from Property (ii), we deduce:

    \vspace{.2cm}

    \noindent{Property {\bf (iii) }} \em For each cell $i$, there exists a first finite time $t_i > 0$ such that $\lim_{t \rightarrow t_i^-}S_i(x_i(t)) = \theta_i $.   Namely, for any $i \in \{1, \ldots, m\}$ there exists a first spiking  instant $t_i > 0$. \em

   \vspace{.3cm}

   Consider the \em minimum natural number $K \geq 1$ \em such that
   \begin{equation}\label{eqn89} K \geq  \frac{  \displaystyle  \max \{   \theta_j \colon \ j \in {\mathcal N}    \} }{\displaystyle \min\{\Delta_{i,j}\colon i,j \in {\mathcal N} \mbox{ such that } i \neq j\}}.\end{equation}
   From Inequalities (\ref{eqnm'>Theorem1}) and (\ref{eqn89}), we deduce:
   \begin{equation} \label{eqn_m_geqK^2} \sqrt{m  } > K.\end{equation}

By hypothesis, $\theta_j >0$ for all $1 \leq j \leq m$. Denote
    \begin{equation}
    \label{eqn88}
    0 < l_j := \frac{\theta_j  }{K} \leq  { \theta_j   } \cdot \frac{\displaystyle \min_{ i \neq j} \Delta_{i,j}}{\displaystyle \max_{1 \leq j \leq m} \ \theta_j    } \leq {\min_{i \neq j  } \Delta_{i,j}}.
    \end{equation}
    For later use, we note the following property:
     \begin{eqnarray}   \label{eqn82aa}\mbox{If } \ 1 \leq h \leq K-1, & \mbox{\ \ \  then \ \ \  } & h \ l_j =   \frac{h }{K} \,  \theta_j     <    \theta_j,  \\ \mbox{ and  if }\ \ \  h= K, & \mbox{\ \ \  then\ \ \  }   &   K \ l_j = \frac{K }{K} \,  \theta_j   =   \theta_j. \label{eqn82a}\end{eqnarray}

   \noindent Assume that at the instant $t > 0$, at leat $K$   cells are spiking, where $K \geq 1$ is the natural number defined by Inequality (\ref{eqn89}). Then, for any other cell $j \in \{1, \ldots, m\}$, applying Equalities (\ref{eqnInteraction}), we have:
   $$S_j(x_j(t)) \geq  \  S_j(x_j(t^-)) \ + \ K \min _{i \neq j } \Delta_{i,j} \geq $$ $$\geq    \  \big(\min _{i \neq j } \Delta_{i,j}\big) \cdot  \frac{\displaystyle \max_{1 \leq j \leq m} \ \theta_j }{\displaystyle \min_{i \neq j } \Delta_{i,j}}  = \  {\max_{1 \leq j \leq m} \ \theta_j \ }     \geq    \theta_j.$$
   Thus, any cell $j$ will also spike at instant $t$, because its satisfaction variable $S_j(x_j(t))$  arrives to     the goal level $\theta_j$. Summarizing, we have proved:

\vspace{.3cm}

   \noindent{Property {\bf (iv) }}   \em If at least $K$   cells   spike  at an instant $t > 0$, then all the cells   spike at the instant $t$. \em

   \vspace{.3cm}

   Now, to end the proof of Assertion (i) it is enough to prove the following:

   \vspace{.3cm}

   \noindent{Assertion {\bf (v) } to be proved: }   \em There exists an instant $T> 0 $ such that   at least $K$   cells spike simultaneously at $T$. \em

\vspace{.3cm}

    Let us take $  t_{1} > 0$ equal to the first positive instant when at least one cell $i_1  $  arrives to its goal level, i.e. \begin{equation}\label{eqn85} S_{i_1 }(x_{i_1 }(t_{1})^-) = \theta_{i_1 }  \ \ \mbox{ for some } \ {i_1 } \in \{1, \ldots, m\}.\end{equation}

    Let us discuss according to two cases: either at least $K$ cells spike at instant $t_{1}$ with the cell $i_1$, or at most $K-1$ do.

\vspace{.2cm}

    \noindent{FIRST CASE: } At least $K$ cells spike at instant $t_{1} > 0$.    Thus Assertion (v) holds.

\vspace{.2cm}

\noindent{SECOND CASE: } At most $K-1$ cells (and at least one cell $i_1$) spike  at instant $t_{1}$. From Inequality (\ref{eqn_m_geqK^2}), there exist at least $m - (K-1) \geq K^2- (K-1) \geq 1 $ cells that do not spike at instant $t_{1}$. Denote by $A_1$   this \em set of non spiking   cells \em at instant $t_{1}$. We have:
$$\#A_1 \geq K^2 - (K-1) \geq 1.$$ Using Inequality (\ref{eqn87}), for each cell $j   \in A_1$  we know that
$S_j(x_j(t_{1}^-)) \geq 0$. Since at instant $t_{1}$ at least the cell $i_1$ spikes, it sends a positive interaction $\Delta_{i_1,j}$ to any cell $j \in A_1$. Applying Inequalities (\ref{eqnInteraction}) and (\ref{eqn88}), we deduce:
$$S_j(x_j(t_{1})) \geq   \ S_j(t_{1}^-) + \Delta_{i_1,j}   \geq     \min_{i \neq j} \Delta_{i,j}    \geq   l_j \ \ \ \ \ \ \forall \ j  \in A_1.$$
In brief:
\begin{equation}\label{eqn83a}S_j(x_j(t_{1})) \geq     l_j  \ \ \ \ \ \ \forall \ j  \in A_1.\end{equation} Denote by $t_{2} > t_{1}$ the first  instant after $t_{1}$ for which at least one cell $i_2$ arrives to its goal level.

Now, we discuss again two cases: either there exist at least $K$ cells that spike at instant $t_{2}$ with the cell $i_2$, or there exists at most  $K-1$ cells that so do. In the first case, Assertion (v) holds. In the second case,
denote by $A_2 \subset A_1$ the set of cells that did not spike at any instant  in  $[0,t_{2}]$.  We   have:
$$\#A_2 \geq \#A_1 - (K-1) \geq K^2 - 2(K-1).$$ Since $t_{2} > t_{1}$, applying Property (ii) and Inequalities (\ref{eqnInteraction}) and  (\ref{eqn83a}), we obtain:
$$S_j(x_j(t_{2}^-)) > S_j(x_j(t_{1})) \geq   l_j \ \ \forall \ j \in A_2.$$
Since at least one cell $i_2$ spikes at instant $t_{2}$, it adds a positive jump $\Delta_{i_2,j}$ to $S_j(x_j(t_{2}^-))$ for all $j \in A_2$. Thus, using Inequality (\ref{eqn88}), we deduce:
\begin{equation}\label{eqn90} \begin{split}S_j(x_j(t_{2})) \geq \  S_j(x_j(t_{2}^-)) + \Delta_{i_2,j}    \geq    l_j + \min_{i \neq j} \Delta_{i,j}   \geq      2l_j   \ \ \ \ \forall \ j \in A_2.\end{split} \end{equation}

 By induction on $h \in \mathbb{N}^+$, if $2 \leq h \leq K$, assume that     $t_{h}$ is the $h$-th. consecutive instant $t_{h} > t_{{h-1}}$ such that at least one cell $i_h$ spikes, and no more than $K-1$ cells have spiked simultaneously at each instant  $ t_{1} < t_{2} \ldots < t_{ {h-1}}$. Denote by $A_h$ the set of cells that have not spiked at those instants and also do not spike  at   instant $t_{h}$.

 We have \begin{equation} \label{equation}\#A_h \geq m - h(K-1) \geq K^2 - h(K-1).\end{equation}
Arguing  by induction as in Inequality (\ref{eqn90}), we obtain:
\begin{equation}
\label{eqn91} S _j(x_j(t_{h})) \geq     h \, l_j  \ \ \  \ \ \forall \ j \in A_h.
\end{equation}
For $h= K$, joining Inequality (\ref{eqn91}) and Equality (\ref{eqn82a}), we deduce:
$$S_j(x_j(t_{K})) \geq  \theta_j  \ \ \forall \ j \in A_K.$$
In other words, we have proved that, if at each instant  $  t_{1} < t_2 < \ldots < t_{{K-1}}$, not more than $K-1$ cells spike simultaneously, then, at instant $t_{K}$ the value of $S_j$ arrives to the goal level $\theta_j$, for any cell $j \in A_K$. This implies that all the cells of $A_K$ spike simultaneously at   instant  $t_{K} >0$. Besides, from Inequality (\ref{equation}) there exist at least $K^2 - K(K-1) = K$ cells in the set $A_K$. So, we have proved that at least $K$ cells spike simultaneously at some instant $t_1 >0$ or $t_2 >t_1$ or $\ldots t_{K-1}$ or $t_K$.   This ends the proof of Assertion (v), as wanted, and thus the proof of Part (a) of Theorem \ref{theoremFullCooperativism} is complete.
\hfill $\Box$

\noindent{\em Proof: }  {\bf Part b) }

Since the hypothesis (\ref{eqnm'>Theorem1}) is a strict inequality, it establishes  an open condition in the space of parameters of the network, endowed with the topology of Definition \ref{definitionGenericNetworks}. Therefore, the  Inequality (\ref{eqnm'>Theorem1}) joint with any dynamical property that is deduced from it, is a robust phenomenon. In particular, due to Part (a) of Theorem \ref{theoremFullCooperativism},  Inequality (\ref{eqnm'>Theorem1}) joint with the eventually periodic synchronization of the spikes,   is a robust phenomenon. This proves Assertion (b) of Theorem \ref{theoremFullCooperativism}. \hfill $\Box$

\subsection
{\bf Proof of  Theorem \ref{TheoremPart(c)}}
\label{proofPart(c)Theorem1}

\noindent{\em Proof: }  From the proof of Part (a) of Theorem \ref{theoremFullCooperativism}, the waiting time $T > 0$ until the spike-synchroniz\-ation  of all the cells occurs, equals the time   that takes the latest  cell, say $i$, to arrive to its goal level $\theta_i$ from its initial state    $ x_i(0)$. The worst case occurs when the initial state $x_i(0)$ of this slowest cell $i$ is such that $S_i(x_i(0))$ takes its lowest possible value 0 (cf. Inequality (\ref{eqn87})). Thus, to consider the worst case, we assume $$S_i(x_i(0)) = 0.$$ Therefore $T$ is smaller or equal than the time constant $t_i$ of the differential equation (\ref{eqnSelfSatisfactionDiffEq}), for the solution $S_i(x_i(t))$ with initial state $S_i(x_i(0))= 0$, such that $S_i(x_i(t_i^-)) = \theta_i$. This is because  during the time-interval $[0, t_i)$, some other cells $j \neq i$ might have spiked. So they might have injected non negative jumps $\Delta_{j,i}$ to the instantaneous value of the variable $S_i$ of the cell $i$.

Then, the worst case if when all those jumps are zero. Summarizing, we have $$T \leq t_i$$ and the worst case occurs if    $S_i$ is only governed by the differential equation (\ref{eqnSelfSatisfactionDiffEq}) for all $t \in [0, t_i)$. Due to the mean value theorem of the differential calculus, there exists a time $\tau_i \in [0, t_i)$ such that:
\begin{equation}\label{eqn51bb}\left . \frac{dS_i}{dt}\right |_{\displaystyle t = \tau_i} = \frac{S_i(x_i(t_i^-)) - S_i(x_i(0))}{t_i} = \frac{\theta_i  }{t_i}.\end{equation}
Using the differential equation (\ref{eqnSelfSatisfactionDiffEq}), we have:
\begin{equation}\label{eqn51cc}\left. \frac{dS_i}{dt}\right |_{\displaystyle t = \tau_i} = g_i  (x_i(\tau)) \geq  \min \{g_i (x_i) \colon x_i \in X_i\}.\end{equation}
Joining Equality (\ref{eqn51bb}) and Inequality (\ref{eqn51cc}), we obtain:
$$T \leq t_i \leq \frac{\theta_i}{\min \{g_i (x_i) \colon x_i \in X_i \} } \leq \max_{i \in {\mathcal N}} \Big \{  \frac{\theta_i}{\min \{g_i (x_i) \colon x_i \in X_i \} } \Big\},$$ proving the first assertion   of Theorem \ref{TheoremPart(c)}.

\vspace{.2cm}

Now, let us prove the second assertion of   Theorem \ref{TheoremPart(c)}, finding an upper bound for the natural spiking period $p$. We revisit the proof of Part (a) of Theorem \ref{theoremFullCooperativism}. We have defined the constant $K$ as the minimum positive natural number that satisfies Inequality (\ref{eqn89}). Thus:
$$K \leq 1 + \frac{  \displaystyle  \max \{   \theta_j \colon \ j \in {\mathcal N}\}    }{\displaystyle \min \{ \Delta_{i,j}\colon i,j \in {\mathcal N} , i \neq j  \}}.$$
On the one hand, by Property (iv) in the proof of Part (a) of Theorem \ref{theoremFullCooperativism}, if $K$ cells spike simultaneously at instant $t_i$ then all the cells spike simultaneously at instant $t_i$. On the other hand, at the end of the proof of Assertion (v), (second case), we found that  from any initial state,   after at most $K$ spikes of some cells (i.e.   at instant $t_K$ as latest), there exist $K$ cells that spike simultaneously. We conclude that, once all the cells had   simultaneously spiked at instant $t_0$, the minimum next instant $t_p >t_0$ for which all  the cells  spike again simultaneously, is such that $p \leq K$. Therefore:
$$p \leq K \leq 1 + \frac{  \displaystyle  \max \{   \theta_j \colon \ j \in {\mathcal N}\}    }{\displaystyle \min \{ \Delta_{i,j}\colon i,j \in {\mathcal N} , i \neq j  \}},$$
ending the proof of   Theorem \ref{TheoremPart(c)}.
\hfill $\Box$

\section{The proof of the amount of information}  \label{sectionEntropy} \label{proofPart(e)Theorem1}

In this Section,  we prove Theorem \ref{TheoremPart(e)}. We compute the total amount of information of   fully cooperative networks with a large number $m$ of cells.

%%%%%%%%%%%%%%%%%%%%%55

\noindent{\em Proof: }  {\em of Theorem } \ref{TheoremPart(e)}

In Part (B) of Lemma \ref{propositionSynchronizationPeriodic}  we have proved that
$$\max_{k \geq 0} \{\#{\mathcal P}_k\} = p,$$
where $p$ is the natural spiking period, i.e the period of the sequence $\{I_n\}_{n \in \mathbb{N}}$ of spiking-codes. Therefore, from Formula (\ref{eqnDefinitionEntropia}), the total amount of information $H$ that the network can generate or process is
$$H= \sup_{k \geq 0} \{\log_2(\#{\mathcal P}_k) \}= \log_2 (\sup_{k \geq 0} \{\#{\mathcal P}_k\}) = \log_2 (\max_{k \geq 0} \{\#{\mathcal P}_k\}) = \log _2 p.$$
(The first equality in the above chain holds because the real function $\log_2(x)$ is increasing on $x \in \mathbb{R}^+$.) We have proved that   $$H= \log_2 p ,$$ which is the first assertion of Theorem \ref{TheoremPart(e)}. Now, let us prove the second assertion. We use the upper bound of the natural spiking period $p$ that was proved in Theorem \ref{TheoremPart(e)}, Formula (\ref{eqnBoundPeriodp}). We conclude that
$$H = \log_2 p \leq \log_2 \Big(1 + \frac{  \displaystyle  \max_{1 \leq j \leq m}   \theta_j     }{\displaystyle \min_{i \neq j } \Delta_{i,j}}\Big),$$
ending the proof of Theorem \ref{TheoremPart(e)}.
\hfill $\Box$

%%%%%%%%%%%%%%%%%%%%%%%%%%%%%
%%%%%%%%%%%%%%%%%%%%%%%%%%%%%%%%%%%%%%%%%%%%%%%%%%%%%%%%%%

\section{The proofs of results on the protection factor}  \label{sectionProtectionFactor}

In this section we prove Proposition \ref{propositionProtectionFactor} and  Theorem \ref{TheoremPart(d)}.

\subsection
 {\bf Proof of Proposition \ref{propositionProtectionFactor}} \label{subsectionProofPrpositionProtectionFactor}

\noindent{\em Proof: }
From Equalities (\ref{eqnInstrinsicRisk}) and (\ref{eqnProtectionFactor}), we obtain
$$\big(1  - P_i \big) \ R_i = \frac{\Big(\displaystyle \theta_i - \Delta_{i,{\mathcal N}}^{(h)}\Big) \   \theta_i}{\ \theta_i\ \ \Big(\displaystyle \max_{1 \leq j \leq m} \theta_j\Big)} = \frac{ \displaystyle    \     \theta_i -\Delta_{i,{\mathcal N}}^{(h)} \  }{  \displaystyle \max_{1 \leq j \leq m} \theta_j }$$
From Equality (\ref{eqnNetRisk}),
$\displaystyle R'_i = \max\Big \{0, \  \min\Big\{1, \   \frac{\displaystyle  \    \theta_i - \Delta_{i, {\mathcal N}}^{(h)} \ }{\displaystyle   \max_{1 \leq j \leq m} \theta_j} \Big\}, $

\noindent and thus: $\displaystyle  \ R'_i = \max \Big \{0, \ \min\big\{1, \ \ \big( 1 - P_i\big) \ R_i \big \} \Big\}.$
\hfill $\Box$

\subsection
 {\bf Proof of Theorem \ref{TheoremPart(d)}} \label{proofPart(d)Theorem1}

\noindent{\em Proof: }
  By hypothesis the network is fully cooperative. Thus, $\Delta_{j,i} >0$ for all $j \neq i$. By Equalities (\ref{eqnNetDelta}) and  (\ref{eqnProtectionFactor}), the protection factor $P_i$ of the network to each cell is positive.

Since $P_i >0,$   we obtain
$\displaystyle  1-P_i  < 1. $
Applying Formula (\ref{eqnFormulaProtectionFactor}) of Lemma $\ref{propositionProtectionFactor}$, and recalling that $0 <R_i \leq 1$ -cf. Inequality (\ref{eqnR_i>0})- we obtain:
 $$R'_i < \max\{0, \min\{1, R_i\}\} = R_i.$$
 Therefore the net risk of death $R'_i$ of the cell $i$ when interacting in the network is strictly smaller than the intrinsic risk $R_i$ that $i$ would have if it were not connected to the network. This ends the proof of Theorem \ref{TheoremPart(d)}.
\hfill $\Box$

\section{Proof of the results on networks with similar cells}
 In this section we prove Theorem \ref{TheoremParts(fgh)}:

\label{proofPart(f)Theorem1}

\noindent{\em Proof: }  {\bf Part a) of Theorem \ref{TheoremParts(fgh)}}

From Part (a) of Theorem \ref{theoremFullCooperativism}, there exists a first instant $t_0 >0$ such that all the cells spike simultaneously at $t_0$. So, to proof Part (a) of Theorem \ref{TheoremParts(fgh)} it is enough to show the following assertion, under the additional hypothesis of Inequality (\ref{eqnhypothesisnueva}):

\vspace{.2cm}

\noindent{Assertion {\bf (vi) } to be proved. }  \em If at some instant $t_0>0$ all the cells spike simultaneously, and if $t_1 > t_0$ is the first   instant after $t_0$ when at least one cell spikes, then all the cells spike simultaneously at $t_1$. \em

\vspace{.2cm}

Fix $i$, one of the cells that spike  at instant $t_1 > t_0 $. By hypothesis, the cell $i$ has  also spiked   at instant $t_0$, but not during the inter-spike-interval $(t_0, t_1)$. Then, applying the reset rule (\ref{eqnSelfSatisfactionReset}), we have:   $$S_i(x_i(t_0)) = 0$$
Since $i$ also spikes at instant $t_1 >t_0$, we have: $$S_i(x_i(t_1^-)) = \theta_i.$$
Therefore:
$$S_i(x_i(t_1^-)) - S_i(x_i(t_0)) = \theta_i.$$
No cell spikes during the time-interval $(t_0, t_1)$. Thus, $S_i$ is governed by the differential equation (\ref{eqnSelfSatisfactionDiffEq}) during such a time-interval. Applying the mean value theorem of the differential calculus, there exists $\tau_i \in (t_0, t_1)$ such that:
\begin{equation} \label{eqn70}\left . \frac {dS_i}{dt} \right|_{\displaystyle t= \tau_i} = \frac{S_i(x_i(t_1^-)) - S_i(x_i(t_0))}{t_1- t_0} = \frac{\theta_i}{t_1- t_0}.\end{equation}
 Using the differential equation   (\ref{eqnSelfSatisfactionDiffEq}), and recalling that $x_i \in X_i$, $X_i$ is compact and $g$ is continuous,  we obtain:
\begin{equation} \label{eqn71}\frac{dS_i}{dt}   = g_i(x_i)  \leq \max\{ g_i (x_i) \colon x_i \in X_i\}. \end{equation}
Joining Equality (\ref{eqn70}) and Inequality (\ref{eqn71}), we deduce:
\begin{equation} \label{eqn73}t_1 - t_0 \geq \frac{\theta_i}{\max\{ g_i (x_i) \colon x_i \in X_i\}} \geq    \frac{\displaystyle  \min_{1 \leq i \leq m} {\theta_i}}{\displaystyle  \max_{1 \leq i \leq m} \max\{ g_i (x_i) \colon x_i \in X_i\}}.\end{equation}

Now denote by $j$ any cell. It spikes at instant $t_0$. Then,  $$S_j(x_j(t_0)) = 0, \ \ 0 \leq S_j(x_j(t_1^-)) \leq \theta_j.   $$

No cell spikes during the time-interval $(t_0, t_i)$. Thus, $S_j$ is governed by the differential equation (\ref{eqnSelfSatisfactionDiffEq}) during such an interval of time. Applying the mean value theorem of the differential calculus, there exists $\tau_j \in (t_0, t_1)$ such that:
\begin{equation} \label{eqn70bb}\left . \frac {dS_j}{dt} \right|_{\displaystyle t= \tau_j} =\frac{S_j(x_j(t_1)^-) - S_j(x_j(t_0))}{t_1 - t_0} =  \frac{S_j(x_j(t_1)^-) }{t_1 - t_0}.\end{equation}

From the differential equation (\ref{eqnSelfSatisfactionDiffEq}),    we deduce:
\begin{equation} \label{eqn71bb}\frac{dS_j}{dt}   = g_j (x_j )\geq  \min  \{ g_j(x_j) : \ \  0 \leq x_j \leq \theta_j  \} \geq \min_{1 \leq i \leq m}\min  \{ g_i(x_i) : \ \  x_i \in X_i  \} . \end{equation}

Joining Inequalities (\ref{eqn70bb}) and  (\ref{eqn71bb}), we deduce:
\begin{equation} \label{eqn73b}t_1- t_0 \leq \frac{S_j(x_j( t_1^-)) }{\displaystyle \min_{1 \leq i \leq m}\min  \{ g_i(x_i) : \ \  x_i \in X_i  \}  \}}. \end{equation}

Joining Inequalities (\ref{eqn73}) and (\ref{eqn73b}), we obtain:
$$
{S_j(x_j( t_1^-)) } \geq    \frac{\displaystyle \Big(\min_{1 \leq i \leq m}\min  \{ g_i(x_i) : \ \  x_i \in X_i  \}  \}\Big) \ \Big({\displaystyle  \min_{1 \leq i \leq m} {\theta_i}} \Big)}{ {\displaystyle  \max_{1 \leq i \leq m} \max\{ g_i (x_i) \colon x_i \in X_i\}}}
 $$
from where

 ${S_j(x_j( t_1^-)) } - \theta_j \geq $

  \hfill $\geq \frac{\displaystyle \Big(\min_{1 \leq i \leq m}\min  \{ g_i(x_i) : \ \  x_i \in X_i  \}  \}\Big) \ \Big({\displaystyle  \min_{1 \leq i \leq m} {\theta_i}} \Big)}{  \displaystyle  \Big(\max_{1 \leq i \leq m} \max\{ g_i (x_i) \colon x_i \in X_i\} \Big) \ \Big({\displaystyle  \max_{1 \leq i \leq m} {\theta_i}} \Big)} \ \cdot \Big({\displaystyle  \max_{1 \leq i \leq m} {\theta_i}} \Big) \ - \ \theta_j \geq $  \begin{equation} \label{eqnCCC}\geq  {\displaystyle  \max_{1 \leq i \leq m} {\theta_i}}  \ \cdot \  \left (\frac{\displaystyle \Big(\min_{1 \leq i \leq m}\min  \{ g_i(x_i) : \ \  x_i \in X_i  \}  \}\Big) \ \Big({\displaystyle  \min_{1 \leq i \leq m} {\theta_i}} \Big)}{  \displaystyle  \Big(\max_{1 \leq i \leq m} \max\{ g_i (x_i) \colon x_i \in X_i\} \Big) \ \Big({\displaystyle  \max_{1 \leq i \leq m} {\theta_i}} \Big)} - 1   \right) \end{equation}
 By hypothesis, Inequality (\ref{eqnhypothesisnueva}) holds. So, the factor at right (between large parenthesis) in Inequality (\ref{eqnCCC}) is bounded from below by $- (\min_{i \neq j} \Delta_{i,j}) / \max_{1 \leq i \leq m} \theta_{i}$. We deduce:
 $${S_j(x_j( t_1^-)) } - \theta_j \geq   - \min_{i \neq j} \Delta_{i,j},
 $$
from where
\begin{equation}
 \label{eqnCCD}S_j(x_j(t_1^-)) + \min_{i \neq j} \Delta_{i,j} \geq \theta_j. \end{equation}

Since at instant $t_1  $ the cell $i$ spikes,   it sends an action $\Delta_{i,j}$ to each cell $j$. From Inequality (\ref{eqnCCD}), and from the interaction rule in Equalities (\ref{eqnInteraction}), we get:
$$S_j(x_j(t_1)) =  S_j(x_j(t_1^-)) + \Delta_{i,j} \geq  S_j(x_j(t_1^-)) + \min_{i \neq j} \Delta_{i,j} \geq \theta_j. $$
Therefore, the variable $S_j$ of any cell $j$ arrives to its  goal level $\theta_j $ at the instant $t_1$. Thus, any cell $j$ spikes at $t_1$, and   Assertion (vi) is proved. This ends the proof of Part (a) of Theorem \ref{TheoremParts(fgh)}.
\hfill $\Box$

\noindent{\em Proof: }    {\bf Part b) of Theorem \ref{TheoremParts(fgh)}}

First, let us prove that the protection factor $P_i$ of the network to each cell $i$, according to Definition \ref{definitionProtectionFactor}, satisfies \begin{equation}\label{eqn44}\min\{1, \ P_i \} = 1   \ \ \ \mbox{ (to be proved). } \end{equation}
In fact, by Formula (\ref{eqnNetDelta}) the net sum of the actions $\Delta_{i, {\mathcal N}}^{(h)}$ that the cell $i$ receives from the other cells of the network during its $h$-th. inter-spike-interval $(t_i^{(h)}, t_i^{(h+1)}]$ is
\begin{equation}\label{eqn45}\Delta_{i, {\mathcal N}}^{(h)} = \sum_{\displaystyle j \neq i, \ \ j \in I_n, \ t_n \in (t_i^{(h)}, t_i^{(h+1)}]} \Delta_{j,i}\end{equation}
By Part (a) of Theorem \ref{TheoremParts(fgh)}, all the cells synchronize spikes after a transitory time $T$. So, for all $h \geq 1$ such that $t_i^{(h+1)} \geq T$, all the cells belong to $I_n$ for the spiking time $t_n = t_i^{(h+1)}$. Thus, in Equality (\ref{eqn45}) the sum at right is realized in $j \in \{1, \ldots, m\}$ such that $j \neq i$. Then,
\begin{equation}\label{eqn46}\Delta_{i, {\mathcal N}}^{(h)} = \sum_{\displaystyle j \neq i} \Delta_{j,i} \geq (m-1) \min_{j \neq i} \Delta_{j,i},\end{equation}
where $m$ is the number of cells of the network.
By hypothesis, Inequality (\ref{eqnm'>Theorem1}) holds. Therefore $m \geq 3$, which implies $m-1 \geq \sqrt m$. We deduce that:
$$m-1 \geq \frac{\displaystyle \max_{1 \leq i \leq m} \theta_i }{\displaystyle \min_{j \neq i} \Delta_{j,i}}.$$
Substituting   in (\ref{eqn46}), we obtain:
$\displaystyle \Delta_{i, {\mathcal N}}^{(h)}   \geq \max_{1 \leq i \leq m} \theta_i  \geq \theta_i $. Thus, $\displaystyle \frac{ \Delta_{i, {\mathcal N}}^{(h)}}{ \theta_i}   \geq   1$.
Now, we apply Formula (\ref{eqnProtectionFactor}):
$\displaystyle P_i = \frac{  \ \Delta_{i, {\mathcal N}}^{(h)} \ }{\theta _i} \geq 1,$
and thus Equality (\ref{eqn44}) is proved.

Second, we apply Lemma \ref{propositionProtectionFactor}. From Equalities  (\ref{eqnFormulaProtectionFactor}) and (\ref{eqn44}), and Inequality (\ref{eqnR_i>0}), we deduce:
$$R'_i = \max\big\{0, \ \min\{1 , \  \ (100- P_i)  R_i\ \}\big \} = \max\{0, \  (1- P_i)  R_i\} = 0,$$
ending the proof of Part (b) of Theorem \ref{TheoremParts(fgh)}.
\hfill $\Box$

\noindent{\em Proof: }    {\bf Part c) of Theorem \ref{TheoremParts(fgh)}}

From Part (a) of Theorem \ref{TheoremParts(fgh)}, under the additional hypothesis stated by Inequality (\ref{eqnhypothesisnueva}), the period $p$ of the sequence $\{I_n\}_{n \in \mathbb{N}}$ of spiking-codes is $p= 1$. From Theorem \ref{TheoremPart(e)}, we know that the total amount of information $H$ that the network can generate of process is $\log_2 p$. Therefore,
$H= \log_2 p = \log_2 1 = 0,$ as wanted.
\hfill $\Box$

%%%%%%%%%%%%%%%%%%%%%%%%%%%%%
%%%%%%%%%%%%%%%%%%%%%%%%%%%5

\section {\bf Proof of Corollary \ref{corollary1}} {\bf (Large full cooperative core)} \label{proofCorollary1} \label{sectionProofOfCorollaries}

\noindent{\em Proof: }
  Since by hypothesis, the number $m$ of cells of the core ${\mathcal N}_1$ satisfies Inequality (\ref{eqn8cor1}), we can repeat all the arguments in the proof of Part (a) of Theorem \ref{theoremFullCooperativism} in Subsection \ref{proofPart(a)Theorem1}, by substituting the network ${\mathcal N}$ by the core ${\mathcal N}_1$. So, we deduce that there exists a strictly increasing sequence $\{t_n\}_{n \geq 0}$ of instants $t_n \rightarrow + \infty$, and a natural period $p \geq 1$, such that:

       At least one cell of ${\mathcal N}_1$ spikes at each instant $t_n$, for all $n \geq 0$.

      No cell of ${\mathcal N}_1$ spikes in the open time-intervals $(t_n, t_{n+1})$.

     All the cells of ${\mathcal N}_1$ spike simultaneously at the instant $t_{hp}$ for all $h \geq 0$.

  \vspace{.2cm}
  Therefore, to prove the eventual periodic synchronization of the whole network ${\mathcal N}$, it is enough to prove that   the cells of ${\mathcal N} $ spike altogether at each instant $t$ such that all the cells of the core ${\mathcal N}_1$ spike.

  In fact, we repeat, with slight changes, the proof of Property (iv) in Subsection \ref{proofPart(a)Theorem1}: We define the \em minimum natural number \em $K \geq 1$ such that
  $$K > \frac{\max\{\theta_j \colon \ j \in {\mathcal N}\}}{\min\{\Delta_{i,j}\colon i \in {\mathcal N}_1, \ j \in {\mathcal N} \ , \ i \neq j\}}.$$
  Due to Inequality (\ref{eqn8cor1}), the number $m$ of the core ${\mathcal N}_1$ satisfies
  \begin{equation} \label{eqn1010}\sqrt m > K.\end{equation}
  Now, we repeat the same arguments of the proof in Subsection \ref{proofPart(a)Theorem1}, starting at Equality (\ref{eqn88}) and ending just before Property (iv):  Substituting
  $\min_{i \neq j} \Delta_{i,j}$ by
  \begin{equation} \label{eqn1011}\min \{\Delta_{i,j} \colon  \ i \in {\mathcal N}_1, \ j \in {\mathcal N} \ , i \neq j\},\end{equation} we deduce

  \noindent Property {\bf (iv)'}: \em If at least $K$ cells of the core ${\mathcal N}_1$ spike at instant $t >0$, then all the cells of the network ${\mathcal N}$ spike at instant $t$. \em

  Since at instant  $t_{hp}$ all the $m$ cells of the core ${\mathcal N}_1$ spike, and $m \geq 3$ satisfies Inequality (\ref{eqn1010}), we obtain:
  $$m > \sqrt m > K.$$
  We conclude that   the cells of the network spike altogether at instants $t_{hp}$ for all $h \geq 0$. Thus, the whole network periodically synchronizes spikes  from any initial state (after a finite transitory time).

  \vspace{.2cm}

  Once we have proved the periodic synchronization of the network, the proof  of the other   assertions  under the hypothesis of Corollary \ref{corollary1}, follow the same arguments in the proofs of Sections \ref{sectionSincronizacion}, \ref{sectionEntropy} and \ref{sectionProtectionFactor},  after substituting $\min_{i \neq j} \Delta_{i,j}$ by the expression (\ref{eqn1011}).
\hfill $\Box$
%\subsection
%{\bf Proof of Corollary \ref{corollary2}}{\bf (The small world principle)} \label{proofCorollary2}
%
%\noindent{\em Proof: } 
%
%Under the hypothesis of Corollary \ref{corollary2}, we repeat the same arguments detailed in Subsections \ref{proofPart(a)Theorem1}, \ref{proofPart(b)Theorem1}, \ref{proofPart(c)Theorem1}, \ref{proofPart(f)Theorem1} \ref{proofPart(e)Theorem1}, \ref{proofPart(h)Theorem1}, \ref{proofPart(d)Theorem1} and \ref{proofPart(g)Theorem1},  substituting $\min_{i \neq j} \Delta_{i,j}$ by the expression $$\frac{\min\{\Delta_{i,j} \colon \ i,j \in {\mathcal N}, \ i \neq j \  \mbox{ such that }  \Delta_{i,j} >0  \}}{\mbox{e-diam}({\mathcal N})}.$$
%In particular, the number $K$ defined by Inequality (\ref{eqn89}) must be substituted by the minimal natural number $K \geq 1$ such that
%$$K \geq \frac{\big (\max \{ \theta_j  \colon \ j \in {\mathcal N} \}\big) \ \big(\mbox{e-diam}({\mathcal N}) \big)}{\min\{\Delta_{i,j} \colon \ i,j \in {\mathcal N}, \ i \neq j \  \mbox{ such that }  \Delta_{i,j} >0  \}}.$$
%\hfill $\Box$
%
%

%%%%%%%%%%%%%%%%%%%%%%%%%%
%%%%%%%%%%%%%%%%%%%%%%%%%%%

%%%%%%%%%%%%%%%%%%%%%%
%%%%%%%%%%%%%%%%%%%%%%%%

%\section*{Acknowledgments}

% You may incorporate your references as follows in your main tex file.
% Using BibTex is not recommended but can be handled.

\vspace{1cm}

\end{document}